\documentclass[a4paper,12pt]{article} 
\usepackage{amsmath} 
\usepackage{amsopn} 
\usepackage{amssymb} 
\usepackage[mathscr]{eucal} 
\usepackage{theorem} 
\usepackage{enumerate} 

\theorembodyfont{\itshape} 
\theoremstyle{plain}
  \newtheorem{theorem}{Theorem}[section] 
  \newtheorem{proposition}[theorem]{Proposition} 
   
  \newtheorem{lemma}[theorem]{Lemma} 
  \newtheorem{problem}[theorem]{Problem}

\theorembodyfont{\rmfamily} 
\theoremstyle{plain}
  \newtheorem{remark}{Remark}[section]

\renewcommand{\theequation}%
           {\thesection.\arabic{equation}}

\begin{document} 

\begin{center} 
{\LARGE The topological holonomy group and} 

\vspace{1mm} 

{\LARGE the complexity of horizontality} 

\vspace{6mm} 

{\Large Naoya {\sc Ando} and Anri {\sc Yonezaki}} 
\end{center} 

\vspace{3mm} 

\begin{quote} 
{\footnotesize \it Abstract} \ 
{\footnotesize Based on \cite{AK}, 
we study the complexity of horizontality 
in each twistor space $\hat{E}_{\varepsilon}$ 
associated with an oriented vector bundle $E$ of rank $4$ 
with a positive-definite metric over the $2$-torus $T^2$, 
and obtain classification of the topological holonomy groups in $SO(3)$. 
We observe that there exist many topological holonomy groups in $SO(3)$ 
generated by two finite order elements 
and equipped with noncommutative pairs 
which consist of infinite order elements. 
We find topological holonomy groups 
which are dense in $SO(4)$. } 
\end{quote} 

\vspace{3mm} 

\section{Introduction}\label{sect:intro} 

\setcounter{equation}{0} 

Let $E$ be an oriented vector bundle of rank 4 over $T^2 =S^1 \times S^1$ 
with $S^1 =\mbox{\boldmath{$R$}} /2\pi \mbox{\boldmath{$Z$}}$. 
Let $h$ be a positive-definite metric of the bundle $E$. 
Then $h$ induces a metric $\hat{h}$ of 
the $2$-fold exterior power $\bigwedge^2\!E$ of $E$, 
and $\bigwedge^2\!E$ is decomposed into 
two orientable subbundles $\bigwedge^2_+\!E$, $\bigwedge^2_-\!E$ of rank $3$. 
The twistor spaces $\hat{E}_{\pm}$ associated with $E$ are 
the unit sphere bundles in $\bigwedge^2_{\pm}\!E$ respectively. 
Let $\nabla$ be an $h$-connection of $E$, i.e., 
a connection satisfying $\nabla h=0$. 
Then $\nabla$ induces a connection $\hat{\nabla}$ of $\bigwedge^2\!E$, 
which is an $\hat{h}$-connection. 
In addition, $\hat{\nabla}$ induces connections of $\bigwedge^2_{\pm}\!E$. 
A section $\Omega$ of $\hat{E}_{\varepsilon}$ ($\varepsilon \in \{ +, -\}$) 
corresponds to a complex structure $I_{\Omega}$ of $E$ preserving $h$, 
and $\Omega$ is horizontal with respect to $\hat{\nabla}$ 
if and only if $I_{\Omega}$ is parallel with respect to $\nabla$. 

As in \cite{AK}, 
a polygonal curve $\gamma$ in the $xy$-plane $\mbox{\boldmath{$R$}}^2$ 
is said to be \textit{normal\/} 
if $\gamma$ is oriented and 
if each of the line segments which form $\gamma$ is contained in 
a coordinate curve in $\mbox{\boldmath{$R$}}^2$ 
so that the images of its two end points by the standard projection 
${\rm pr}\,: \mbox{\boldmath{$R$}}^2 \longrightarrow T^2$ coincide with 
each other. 
Let ${\rm NPC}$ denote the set of normal polygonal curves 
in $\mbox{\boldmath{$R$}}^2$ which start from $(0, 0)$. 
For each $\gamma \in {\rm NPC}$ and each element $\omega$ 
of the fiber $\hat{E}_{\varepsilon , {\rm pr} (0, 0)}$ 
of $\hat{E}_{\varepsilon}$ on ${\rm pr} (0, 0)$, 
let $\omega_{\gamma}$ be the element 
of $\hat{E}_{\varepsilon , {\rm pr} (0, 0)}$ 
given by the parallel transport along $\gamma$ 
with respect to $\hat{\nabla}$. 
Then we define the \textit{complexity\/} of the horizontality 
in $\hat{E}_{\varepsilon}$ by 
a subset $X(\omega ):=\{ \omega_{\gamma} \ | \ \gamma \in {\rm NPC} \}$ 
of $\hat{E}_{\varepsilon , {\rm pr}(0,0)}$. 
Therefore the complexity of the horizontality is expressed by 
not only the cardinality of $X(\omega )$ 
but also the placement of the elements 
in $\hat{E}_{\varepsilon , {\rm pr}(0,0)}$. 

By definition, 
$X(\omega )$ is given by 
an at most countable subgroup $G_{\hat{\nabla} , \varepsilon}$ of $SO(3)$. 
In the present paper, 
this subgroup is called the \textit{topological holonomy group\/}  
of $\hat{\nabla}$ in $\bigwedge^2_{\varepsilon}\!E$ at ${\rm pr} (0, 0)$. 
The topological holonomy group $G_{\hat{\nabla} , \varepsilon}$ is 
a subgroup of the holonomy group 
given by fixed two oriented circles which generate 
the fundamental group $\pi_1 (T^2 )$ of $T^2$. 
If $E$ admits a hyperK\"{a}hler structure, 
then one of $\hat{E}_{\pm}$ is a product bundle such that 
each constant section is horizontal, and therefore 
the corresponding topological holonomy group consists of 
only the identity matrix $E_3$. 
If $\hat{E}_{\varepsilon}$ has a horizontal section, 
then $G_{\hat{\nabla} , \varepsilon}$ is isomorphic to 
a finite or countably infinite subgroup of $SO(2)$. 
If $\hat{\nabla}$ is flat, 
then $G_{\hat{\nabla} , \varepsilon}$ is abelian 
and therefore we have a representation of $\pi_1 (T^2 )$ in $SO(3)$. 
Moreover, noticing that at least one eigenvalue of any element of $SO(3)$ 
is equal to $1$, 
we see that there exists an at most double-valued horizontal section 
of $\hat{E}_{\varepsilon}$. 
On the other hand, 
in most cases, $\hat{E}_{\varepsilon}$ does not have 
any horizontal sections. 

We say that the horizontality has 
\textit{finite\/} (respectively, \textit{infinite\/}) \textit{complexity\/} 
for $\omega \in \hat{E}_{\varepsilon , {\rm pr}(0,0)}$ 
if $X(\omega )$ is finite (respectively, infinite). 
If the horizontality has finite complexity for $\omega$, 
then we denote $\sharp X(\omega )$ by $d(\omega )$ or $d$, and 
we say that the horizontality has finite complexity for $\omega$ 
of \textit{degree\/} $d=d(\omega )$. 
Based on \cite{AK}, 
it is observed that 
if the horizontality has finite complexity of degree $d>2$ 
for an element $\omega_0 \in \hat{E}_{\varepsilon , {\rm pr} (0, 0)}$, 
then $X(\omega )$ is finite 
for any $\omega \in \hat{E}_{\varepsilon , {\rm pr} (0, 0)}$, 
and therefore 
the topological holonomy group $G_{\hat{\nabla} , \varepsilon}$ is 
a non-trivial finite subgroup of $SO(3)$. 
Such a subgroup is isomorphic to 
a cyclic group, 
a regular dihedral group, 
the alternating groups of degree $4$, $5$ or 
the symmetric   group  of degree $4$ (\cite{klein}).  
On the other hand, 
if the horizontality has finite complexity of degree $d=1$ or $2$ 
for an element $\omega_0 \in \hat{E}_{\varepsilon , {\rm pr} (0, 0)}$, 
then it is possible that $X(\omega )$ is infinite 
for another element $\omega \in \hat{E}_{\varepsilon , {\rm pr} (0, 0)}$ 
(see Remark~\ref{rem:A}). 

Let $\hat{C}_1$, $\hat{C}_2$ be two elements of $SO(3)$. 
We denote by ${\rm ord}\,(\hat{C}_1 )$, ${\rm ord}\,(\hat{C}_2 )$ 
the orders of $\hat{C}_1$, $\hat{C}_2$ respectively. 
Then ${\rm ord}\,(\hat{C}_k )$ ($k=1, 2$) are valued 
in $\mbox{\boldmath{$N$}} \cup \{ \infty \}$. 
In the present paper, 
we have no need to discuss 
the case where one of given two elements of $SO(3)$ is $E_3$. 
Therefore in almost all situations, 
we suppose $\infty \geq {\rm ord}\,(\hat{C}_1 ) 
                   \geq {\rm ord}\,(\hat{C}_2 ) \geq 2$. 
According to Selberg's lemma (\cite[Theorem 4.1 of Chapter 5]{cassels}, 
\cite{selberg}), 
a finitely generated and infinite subgroup of $GL(n, \mbox{\boldmath{$R$}})$ 
has an element of infinite order. 
In particular, 
since the topological holonomy group $G_{\hat{\nabla} , \varepsilon}$ is 
finitely generated, 
if   $G_{\hat{\nabla} , \varepsilon}$ is infinite, 
then $G_{\hat{\nabla} , \varepsilon}$ has an element of infinite order. 
Suppose that $G_{\hat{\nabla} , \varepsilon}$ has 
a noncommutative pair $(\hat{C}_1 , \hat{C}_2 )$ 
with $\infty ={\rm ord}\,(\hat{C}_1 ) 
      \geq    {\rm ord}\,(\hat{C}_2 ) \geq 3$. 
Then we will prove that $X(\omega )$ is dense 
in $\hat{E}_{\varepsilon , {\rm pr}(0,0)}$ 
for each $\omega \in \hat{E}_{\varepsilon , {\rm pr}(0,0)}$ 
(Theorem~\ref{thm:main}). 
Moreover, based on Theorem~\ref{thm:main}, 
we will prove that for such an $h$-connection $\nabla$, 
$G_{\hat{\nabla} , \varepsilon}$ is 
a dense subset of $SO(3)$ (Theorem~\ref{thm:main2}). 
Notice that a condition ${\rm ord}\,(\hat{C}_2 ) \geq 3$ can not be deleted 
for Theorems~\ref{thm:main} and \ref{thm:main2} 
(Propositions~\ref{pro:pipi/2}, \ref{pro:pipi}), 
while there exist topological holonomy groups 
in $\bigwedge^2_{\varepsilon}\!E$ 
generated by an infinite order element and an order two element, 
and equipped with noncommutative pairs 
satisfying $\infty ={\rm ord}\,(\hat{C}_1 ) 
            \geq    {\rm ord}\,(\hat{C}_2 ) \geq 3$  
(Proposition~\ref{pro:2irr}). 
In Section~\ref{sect:ord>2thg}, 
noticing the minimal polynomials of algebraic numbers and 
cyclotomic polynomials (see \cite[Chapter 3]{nagata}), 
we will find many examples of 
topological holonomy groups in $\bigwedge^2_{\varepsilon}\!E$ 
generated by two finite order elements 
and equipped with noncommutative pairs 
which consist of infinite order elements 
(see Proposition~\ref{pro:inftyinfty}, Remark~\ref{rem:case(i)infty}). 
In addition, we will observe that 
there exist uncountably many such topological holonomy groups 
in $\bigwedge^2_{\varepsilon}\!E$ (Theorem~\ref{thm:inftyinfty}). 
We will obtain classification of all the topological holonomy groups 
in $\bigwedge^2_{\varepsilon}\!E$ (Theorem~\ref{thm:classification}). 
This classification is supported by 
Propositions~\ref{pro:pipi/2}, \ref{pro:pipi}, \ref{pro:2irr}, 
Theorem~\ref{thm:main2} and Selberg's lemma, 
and it can be understood 
in terms of the Tits alternative given in \cite{tits} 
(see Remark~\ref{rem:tits}). 

We will define the topological holonomy group $G_{\nabla}$ 
of $\nabla$ at ${\rm pr} (0, 0)$. 
This is an at most countable subgroup of $SO(4)$, and 
closely related to $G_{\hat{\nabla} , \pm}$ 
by the double covering $SO(4)\longrightarrow SO(3)\times SO(3)$. 
In particular, 
$G_{\nabla}$ is finite if and only if 
both of $G_{\hat{\nabla} , \pm}$ are finite 
(Proposition~\ref{pro:finite}). 
Moreover, 
referring to the proofs of Theorems~\ref{thm:main}, \ref{thm:main2}, 
we will find examples of topological holonomy groups 
which are dense in $SO(4)$ (Theorem~\ref{thm:main3}). 

We will define the topological holonomy group $G_{\nabla}$ 
of an $h$-connection $\nabla$ of a Hermitian vector bundle $(E, h)$ 
over $T^2$ of complex rank $2$. 
Then $G_{\nabla}$ is an at most countable subgroup of $U(2)$. 
Suppose that $G_{\nabla}$ is contained in $SU(2)$. 
By the double covering $\Phi :SU(2)\longrightarrow SO(3)$, 
$G_{\nabla}$ is finite if and only if $\Phi (G_{\nabla} )$ is finite 
(Proposition~\ref{pro:finite2}), 
which is an analogue of Proposition~\ref{pro:finite}. 
We can refer to \cite[Chapter 3]{blichfeldt} 
for finite subgroups of $SU(2)$, 
and \cite{hitchin} for constructions of 
asymptotically locally Euclidean spaces 
by considering the quotients of $\mbox{\boldmath{$C$}}^2$ 
by finite subgroups of $SU(2)$. 
Referring to the last part of the proof of Theorem~\ref{thm:main3}, 
we will see that $G_{\nabla}$ is dense in $SU(2)$ 
if and only if $\Phi (G_{\nabla} )$ is dense in $SO(3)$ 
(Proposition~\ref{pro:dense2}). 
If $G_{\nabla}$ has 
a noncommutative pair $(\hat{B}_1 , \hat{B}_2 )$ 
with ${\rm ord}\,(\hat{B}_1 )=\infty$ 
and  ${\rm ord}\,(\hat{B}_2 )\not= 2, 4$, 
then it will be shown that $G_{\nabla}$ is dense in $SU(2)$ 
(Theorem~\ref{thm:main4}), 
which is an analogous result of Theorem~\ref{thm:main2} for $G_{\nabla}$. 
Finally, referring to the proof of Theorem~\ref{thm:main3}, 
we will prove that there exist $h$-connections of $(E, h)$ such that 
the topological holonomy groups are dense in $U(2)$ 
(Theorem~\ref{thm:main5}). 

\begin{remark} 
Related to horizontal sections, 
it is well-known that 
the complex structure of a K\"{a}hler surface 
corresponds to a horizontal section of a twistor space 
associated with the tangent bundle, and 
one of the twistor lifts of an isotropic minimal surface 
in an oriented Riemannian $4$-manifold is 
a horizontal section of a twistor space 
associated with the pull-back bundle 
(see \cite{bryant}, \cite{ES}, \cite{friedrich} 
for isotropic minimal surfaces). 
\end{remark} 

\begin{remark} 
See \cite{BG}, \cite{BG2}, \cite{kuranishi2} 
for existence and construction of dense subgroups in Lie groups. 
\end{remark} 

\section{Oriented metric vector bundles 
of rank 4 over tori}\label{sect:omvboverT2} 

\setcounter{equation}{0} 

Let $E$, $h$ be as in the beginning of Section~\ref{sect:intro}. 
Let ${\rm pr}\,: \mbox{\boldmath{$R$}}^2 \longrightarrow T^2$ be 
the standard projection as in the second paragraph 
of Section~\ref{sect:intro}. 
Then $h$ induces a metric of the pull-back bundle ${\rm pr}^*\!E$ 
of $E$ by ${\rm pr}$, 
which is also denoted by $h$. 
Let $\Tilde{v}_1$, $\Tilde{v}_2$, $\Tilde{v}_3$, $\Tilde{v}_4$ form 
an orthonormal basis of ${\rm pr}^*\!E_{(0, 0)}$ with respect to $h$ 
such that $(\Tilde{v}_1 , \Tilde{v}_2 , \Tilde{v}_3 , \Tilde{v}_4 )$ gives 
the orientation. 
Let $\nabla$ be an $h$-connection of $E$. 
Then $\nabla$ induces a connection of ${\rm pr}^*\!E$, 
also denoted by $\nabla$. 
For $i=1, 2, 3, 4$, 
let $\Tilde{\xi}_{i, x}$, $\Tilde{\xi}_{i, y}$ be 
parallel sections with respect to $\nabla$ of 
the restrictions ${\rm pr}^*\!E |_{l_x}$, ${\rm pr}^*\!E |_{l_y}$ 
of ${\rm pr}^*\!E$ on 
\begin{equation*} 
l_x =\{ (x, 0) \ | \ x\in \mbox{\boldmath{$R$}} \} , \quad 
l_y =\{ (0, y) \ | \ y\in \mbox{\boldmath{$R$}} \} 
\end{equation*} 
respectively satisfying $\Tilde{\xi}_{i, x} (0, 0) 
                        =\Tilde{\xi}_{i, y} (0, 0) 
                        =\Tilde{v}_i$. 
Then identifying ${\rm pr}^*\!E_{(x+2\pi m, y+2\pi n)}$ 
with ${\rm pr}^*\!E_{(x, y)}$ 
for $m$, $n\in \mbox{\boldmath{$Z$}}$, 
we see that there exist elements $A_1$, $A_2 \in SO(4)$ satisfying 
\begin{equation} 
\begin{split} 
&  (\Tilde{\xi}_{1, x} (a+2\pi , 0)  \ 
    \Tilde{\xi}_{2, x} (a+2\pi , 0)  \ 
    \Tilde{\xi}_{3, x} (a+2\pi , 0)  \ 
    \Tilde{\xi}_{4, x} (a+2\pi , 0)) \\ 
& =(\Tilde{\xi}_{1, x} (a,       0)  \ 
    \Tilde{\xi}_{2, x} (a,       0)  \ 
    \Tilde{\xi}_{3, x} (a,       0)  \ 
    \Tilde{\xi}_{4, x} (a,       0))
     A_1 ,                         \\ 
&  (\Tilde{\xi}_{1, y} (0, b+2\pi )  \ 
    \Tilde{\xi}_{2, y} (0, b+2\pi )  \ 
    \Tilde{\xi}_{3, y} (0, b+2\pi )  \ 
    \Tilde{\xi}_{4, y} (0, b+2\pi )) \\ 
& =(\Tilde{\xi}_{1, y} (0, b      )  \ 
    \Tilde{\xi}_{2, y} (0, b      )  \ 
    \Tilde{\xi}_{3, y} (0, b      )  \ 
    \Tilde{\xi}_{4, y} (0, b      )) 
     A_2 
\end{split} 
\label{AxAy} 
\end{equation} 
for any $a$, $b\in \mbox{\boldmath{$R$}}$. 
The \textit{topological holonomy group\/} of $\nabla$ 
at ${\rm pr} (0, 0)$ is the subgroup $G_{\nabla}$ of $SO(4)$ 
generated by $A_1$, $A_2$. 
Then $G_{\nabla}$ is uniquely determined 
up to a conjugate subgroup of $SO(4)$. 

We set $\overline{E} :=T^2 \times \mbox{\boldmath{$R$}}^4$. 
This is a product bundle over $T^2$. 
The natural inner product of $\mbox{\boldmath{$R$}}^4$ gives 
a metric $h$ of the bundle $\overline{E}$. 
We will prove 

\begin{proposition}\label{pro:prod}
For arbitrarily given two elements $A_1$, $A_2$ of $SO(4)$, 
there exists an $h$-connection $\nabla$ of $\overline{E}$ 
satisfying \eqref{AxAy}. 
\end{proposition} 

\vspace{3mm} 

\par\noindent 
\textit{Proof} \ 
Let $\exp$ denote the exponential map for $SO(4)$. 
Then $\exp$ is a surjective map from the Lie algebra of $SO(4)$ 
onto $SO(4)$. 
Therefore for $A_1$, $A_2 \in SO(4)$, 
there exist elements $P_1$, $P_2$ of the Lie algebra of $SO(4)$ 
satisfying $A_1 =\exp (-2\pi P_1 )$, $A_2 =\exp (-2\pi P_2 )$. 
Let $\xi$ be a section of $\overline{E}$. 
Then $\xi$ is considered to be a section of ${\rm pr}^*\!\overline{E}$. 
In addition, $\xi$ is considered to be 
an $\mbox{\boldmath{$R$}}^4$-valued, doubly periodic function 
on $\mbox{\boldmath{$R$}}^2$. 
We set 
\begin{equation*}
\nabla \xi :=d\xi +P_1 \xi dx +P_2 \xi dy. 
\end{equation*} 
Then we can define by $\nabla$ a connection of $\overline{E}$. 
Since $P_1$, $P_2$ are alternating, 
$\nabla$ is an $h$-connection. 
Let $\Tilde{\xi}$ be a parallel section of the restriction 
of ${\rm pr}^*\!\overline{E}$ on $l_x$. 
Then $\Tilde{\xi}$ is represented as $\Tilde{\xi} =\exp (-xP_1 )\xi_0$ 
on $l_x$ for a constant vector $\xi_0 \in \mbox{\boldmath{$R$}}^4$. 
Therefore we obtain the first relation in \eqref{AxAy}. 
Similarly, we obtain the second relation in \eqref{AxAy}. 
Hence we obtain Proposition~\ref{pro:prod}. 
\hfill 
$\square$ 

\vspace{3mm} 

\begin{remark} 
Let $\mathcal{D}_1$, $\mathcal{D}_2$ be two one-dimensional distributions 
on $T^2$ satisfying 
\begin{itemize} 
\item[{\rm (i)}]{each integral curve of $\mathcal{D}_i$ 
for $i\in \{ 1, 2\}$ is homeomorphic to $S^1$,} 
\item[{\rm (ii)}]{for integral curves $C_1$, $C_2$ 
of $\mathcal{D}_1$, $\mathcal{D}_2$ respectively, 
the intersection $C_1 \cap C_2$ consists of a unique point.} 
\end{itemize} 
Then referring to the above proof of Proposition~\ref{pro:prod}, 
we observe that for the pair $(\mathcal{D}_1 , \mathcal{D}_2 )$ 
and two elements $A_1$, $A_2$ of $SO(4)$, 
there exists an $h$-connection $\nabla$ of $\overline{E}$ 
such that the topological holonomy group $G_{\nabla}$ of $\nabla$ 
at each point of $T^2$ with respect to the two circles $C_1$, $C_2$ 
as in the above (ii) through the point is generated by $A_1$, $A_2$. 
\end{remark} 

\section{The twistor spaces associated with 
oriented metric vector bundles over tori}\label{sect:twspriemvbt} 

\setcounter{equation}{0} 

Let $e_1$, $e_2$, $e_3$, $e_4$ form a local orthonormal frame field 
of $(E, h)$. 
Suppose that $(e_1 , e_2 , e_3 , e_4 )$ gives the orientation of $E$. 
We set 
\begin{equation*} 
\begin{split} 
  \Omega_{\pm , 1} & 
:=\dfrac{1}{\sqrt{2}} (e_1 \wedge e_2 \pm e_3 \wedge e_4 ), \\ 
  \Omega_{\pm , 2} & 
:=\dfrac{1}{\sqrt{2}} (e_1 \wedge e_3 \pm e_4 \wedge e_2 ), \\ 
  \Omega_{\pm , 3} & 
:=\dfrac{1}{\sqrt{2}} (e_1 \wedge e_4 \pm e_2 \wedge e_3 ). 
\end{split} 
\end{equation*} 
Then $\Omega_{\pm , 1}$, $\Omega_{\pm , 2}$, $\Omega_{\pm , 3}$ form 
local orthonormal frame fields of $\bigwedge^2_{\pm}\!E$ respectively 
with respect to the metric $\hat{h}$ 
as in the beginning of Section~\ref{sect:intro}. 
The orientations of $\bigwedge^2_{\pm}\!E$ 
given by $(\Omega_{\pm , 1} , \Omega_{\pm , 2} , \Omega_{\pm , 3} )$ 
do not depend on the choice of $(e_1 , e_2 , e_3 , e_4 )$. 

Let $\Tilde{\omega}_{\pm , 1}$, $\Tilde{\omega}_{\pm , 2}$, 
    $\Tilde{\omega}_{\pm , 3}$ form 
orthonormal bases of $\bigwedge^2_{\pm}\!{\rm pr}^*\!E_{(0, 0)}$ respectively 
such that $(\Tilde{\omega}_{\pm , 1} , \Tilde{\omega}_{\pm , 2} , 
            \Tilde{\omega}_{\pm , 3} )$ give the orientations. 
Then $\Tilde{\omega}_{\pm , 1}$, $\Tilde{\omega}_{\pm , 2}$, 
     $\Tilde{\omega}_{\pm , 3}$ are elements 
of $\widehat{{\rm pr}^*\!E}_{\pm , (0, 0)}$ respectively. 
Let $\hat{\nabla}$ be as in the first paragraph 
of Section~\ref{sect:intro}. 
For $i\in \{ 1, 2, 3\}$, 
let $\Tilde{\Omega}_{\pm , i, x}$, $\Tilde{\Omega}_{\pm , i, y}$ be 
horizontal sections with respect to $\hat{\nabla}$ of 
the restrictions $\widehat{{\rm pr}^*\!E}_{\pm} |_{l_x}$, 
                 $\widehat{{\rm pr}^*\!E}_{\pm} |_{l_y}$ 
of $\widehat{{\rm pr}^*\!E}_{\pm}$ on $l_x$, $l_y$ 
respectively satisfying $\Tilde{\Omega}_{\pm , i, x} (0, 0) 
                        =\Tilde{\Omega}_{\pm , i, y} (0, 0) 
                        =\Tilde{\omega}_{\pm , i}$. 
Identifying $\widehat{{\rm pr}^*\!E}_{\pm , (x+2\pi m, y+2\pi n)}$ 
with $\widehat{{\rm pr}^*\!E}_{\pm , (x, y)}$ 
for $m$, $n\in \mbox{\boldmath{$Z$}}$, 
we see that there exist elements $C_{\pm , 1}$, $C_{\pm , 2} \in SO(3)$ 
satisfying 
\begin{equation} 
\begin{split} 
&  (\Tilde{\Omega}_{\pm , 1, x} (a+2\pi , 0)  \ 
    \Tilde{\Omega}_{\pm , 2, x} (a+2\pi , 0)  \ 
    \Tilde{\Omega}_{\pm , 3, x} (a+2\pi , 0)) \\ 
& =(\Tilde{\Omega}_{\pm , 1, x} (a,       0)  \ 
    \Tilde{\Omega}_{\pm , 2, x} (a,       0)  \ 
    \Tilde{\Omega}_{\pm , 3, x} (a,       0)) 
     C_{\pm , 1} ,                            \\ 
&  (\Tilde{\Omega}_{\pm , 1, y} (0, b+2\pi )  \ 
    \Tilde{\Omega}_{\pm , 2, y} (0, b+2\pi )  \ 
    \Tilde{\Omega}_{\pm , 3, y} (0, b+2\pi )) \\ 
& =(\Tilde{\Omega}_{\pm , 1, y} (0, b      )  \ 
    \Tilde{\Omega}_{\pm , 2, y} (0, b      )  \ 
    \Tilde{\Omega}_{\pm , 3, y} (0, b      )) 
     C_{\pm , 2} 
\end{split} 
\label{CxCy} 
\end{equation} 
for any $a$, $b\in \mbox{\boldmath{$R$}}$. 
Noticing the double covering $SO(4)\longrightarrow SO(3)\times SO(3)$, 
we see that $C_{\pm , k}$ are determined by $A_k$ 
and    that $A_k$ is determined by $C_{\pm , k}$ up to a sign 
for $k\in \{ 1, 2\}$. 
Referring to \cite{AK}, 
we suppose that $\Tilde{\Omega}_{\pm , 1, x}$ are periodic sections 
with period $2\pi$. 
The \textit{topological holonomy group\/} of $\hat{\nabla}$ 
in $\bigwedge^2_{\varepsilon}\!E$ for $\varepsilon \in \{ +, -\}$ 
at ${\rm pr} (0, 0)$ is the subgroup $G_{\hat{\nabla} , \varepsilon}$ 
of $SO(3)$ generated by $C_{\varepsilon , 1}$, $C_{\varepsilon , 2}$. 
Then $G_{\hat{\nabla} , \varepsilon}$ is uniquely determined 
up to a conjugate subgroup of $SO(3)$. 

By Proposition~\ref{pro:prod} together with 
the double covering $SO(4)\longrightarrow SO(3)\times SO(3)$, 
we obtain 

\begin{proposition}\label{pro:prod2}
For arbitrarily given elements $C_{\pm , 1}$, $C_{\pm , 2}$ of $SO(3)$, 
there exists an $h$-connection $\nabla$ 
of $\overline{E} =T^2 \times \mbox{\boldmath{$R$}}^4$ 
satisfying \eqref{CxCy}. 
\end{proposition} 

Let $\omega_{\varepsilon , i}$ ($i=1, 2, 3$) be 
the elements of $\hat{E}_{\varepsilon , {\rm pr}(0, 0)}$ 
given by $\Tilde{\omega}_{\varepsilon , i}$. 
We represent an element $\omega$ of $\hat{E}_{\varepsilon , {\rm pr}(0, 0)}$ 
as $\omega =\sum^3_{i=1} c^i \omega_{\varepsilon , i}$ 
with $\sum^3_{i=1} (c^i )^2 =1$. 
The \textit{complexity\/} of horizontality in $\hat{E}_{\varepsilon}$ 
for $\omega$ is given by 
\begin{equation} 
X(\omega ) 
 =\{ (\omega_{\varepsilon , 1} \ 
      \omega_{\varepsilon , 2} \ 
      \omega_{\varepsilon , 3} )Cp \ | \ 
       C\in G_{\hat{\nabla} , \varepsilon} \} , 
\label{ch} 
\end{equation} 
where $p:=\sum^3_{i=1} c^i p_i$ 
and $p_1 :={}^t [1 \ 0 \ 0]$, 
    $p_2 :={}^t [0 \ 1 \ 0]$, 
    $p_3 :={}^t [0 \ 0 \ 1]$. 
In particular, 
the complexity of 
horizontality for $\omega_{\varepsilon , i}$ ($i\in \{ 1, 2, 3\}$) is 
given by 
\begin{equation*} 
X(\omega_{\varepsilon , i} ) 
 =\{ (\omega_{\varepsilon , 1} \ 
      \omega_{\varepsilon , 2} \ 
      \omega_{\varepsilon , 3} )Cp_i \ | \ 
       C\in G_{\hat{\nabla} , \varepsilon} \} . 
\end{equation*} 
Noticing \eqref{ch}, we set 
\begin{equation} 
X(p):=\{ Cp \ | \ C\in G_{\hat{\nabla} , \varepsilon} \} . 
\label{X(p)} 
\end{equation} 

For a real number $\theta \in \mbox{\boldmath{$R$}}$, 
we set 
\begin{equation} 
C(\theta ):=\left[ \begin{array}{ccc} 
                    1 & 0           &  0           \\ 
                    0 & \cos \theta & -\sin \theta \\ 
                    0 & \sin \theta &  \cos \theta 
                     \end{array} 
            \right] . 
\label{Ctheta} 
\end{equation} 
For each element $C$ of $SO(3)$, 
there exist $\theta \in [0, 2\pi )$ and $V\in SO(3)$ 
satisfying $C=VC(\theta ){}^t V$. 
For $\theta \in [0, 2\pi )$ and $V=[v^i_j ]\in SO(3)$, 
$VC(\theta ){}^t V$ is determined by $\theta$ and the first column of $V$: 
\begin{equation*} 
VC(\theta ){}^t V=\exp (\theta (v^1_1 J_1 +v^2_1 J_2 +v^3_1 J_3 )), 
\end{equation*} 
where $\exp$ is the exponential map for $SO(3)$ 
and $J_1$, $J_2$, $J_3$ are elements of the Lie algebra of $SO(3)$ 
defined by 
\begin{equation} 
J_1   :=\left[ \begin{array}{ccc} 
                0 & 0 &  0 \\ 
                0 & 0 & -1 \\ 
                0 & 1 &  0 
                 \end{array} 
        \right] , \ 
J_2   :=\left[ \begin{array}{ccc} 
                0 & 0 & 1 \\ 
                0 & 0 & 0 \\ 
               -1 & 0 & 0 
                 \end{array} 
        \right] , \ 
J_3   :=\left[ \begin{array}{ccc} 
                0 & -1 & 0 \\ 
                1 &  0 & 0 \\ 
                0 &  0 & 0 
                 \end{array} 
        \right] . 
\label{J} 
\end{equation} 
One of the eigenvalues of $C$ is given by $\lambda =1$. 
The first column ${}^t [v^1_1 \ v^2_1 \ v^3_1 ]$ of $V$ is 
an eigenvector of $C=VC(\theta ){}^t V$ 
corresponding to the eigenvalue $\lambda =1$. 
Since we suppose that $\Tilde{\Omega}_{\pm , 1, x}$ are periodic sections 
with period $2\pi$, 
there exists a number $\theta_{\varepsilon , 1} \in [0, 2\pi )$ 
for each $\varepsilon \in \{ +, -\}$ 
satisfying $C_{\varepsilon , 1} =C(\theta_{\varepsilon , 1} )$. 
Then $A_1$ is given by two elements of $SO(2)$. 
Let $\theta_{\varepsilon , 2} \in [0, 2\pi )$ 
and $V_{\varepsilon , 2} \in SO(3)$ 
satisfy $C_{\varepsilon , 2} 
        =V_{\varepsilon , 2} C(\theta_{\varepsilon , 2} ) 
    {}^t V_{\varepsilon , 2}$. 
Suppose that $C_{\varepsilon , k}$ does not coincide with 
the identity matrix $E_3$ 
for any $\varepsilon \in \{ +, -\}$ and any $k\in \{ 1, 2\}$. 
For $\varepsilon \in \{ +, -\}$, 
let $\phi_{\varepsilon} \in [0, \pi /2]$ be the angle 
between $p_1 :={}^t [1 \ 0 \ 0]$ 
and an eigenvector of $C_{\varepsilon , 2}$ corresponding to $\lambda =1$. 
We can suppose that $\phi_{\varepsilon}$ is the angle 
between $p_1$ and the first column of $V_{\varepsilon , 2}$. 
Then $V_{\varepsilon , 2}$ is represented as 
\begin{equation*} 
\begin{split} 
V_{\varepsilon , 2} 
& =\exp ( \phi_{\varepsilon} (-\sin \gamma_{\varepsilon} J_2 
                              +\cos \gamma_{\varepsilon} J_3 )) \\ 
& =\left[ \begin{array}{ccc} 
           \cos \phi_{\varepsilon} & 
          -\sin \phi_{\varepsilon} \cos \gamma_{\varepsilon} & 
          -\sin \phi_{\varepsilon} \sin \gamma_{\varepsilon} \\ 
           \sin \phi_{\varepsilon} \cos \gamma_{\varepsilon} & 
           \sin^2 \gamma_{\varepsilon} 
          +\cos \phi_{\varepsilon} \cos^2 \gamma_{\varepsilon} &  
          (\cos \phi_{\varepsilon} -1) 
           \cos \gamma_{\varepsilon} \sin \gamma_{\varepsilon} \\ 
           \sin \phi_{\varepsilon} \sin \gamma_{\varepsilon} & 
          (\cos \phi_{\varepsilon} -1) 
           \cos \gamma_{\varepsilon} \sin \gamma_{\varepsilon} & 
           \cos^2 \gamma_{\varepsilon} 
          +\cos \phi_{\varepsilon} \sin^2 \gamma_{\varepsilon}  
            \end{array} 
   \right] 
\end{split} 
\end{equation*} 
for $\gamma_{\varepsilon} \in [0, 2\pi )$, 
and ${}^t [0 \ -\sin \gamma_{\varepsilon} \ \cos \gamma_{\varepsilon} ]$ is 
an eigenvector of $V_{\varepsilon , 2}$ corresponding to $\lambda =1$. 
In the following, we suppose $\phi_{\varepsilon} \in (0, \pi /2]$. 

\section{Finite topological holonomy groups}\label{sect:finitethg} 

\setcounter{equation}{0} 

Suppose that $G_{\hat{\nabla} , \varepsilon}$ is finite. 
Since we suppose that $C_{\varepsilon , k}$ does not coincide with $E_3$ 
and that $\phi_{\varepsilon} \in (0, \pi /2]$, 
$G_{\hat{\nabla} , \varepsilon}$ is isomorphic to 
one of the following (\cite{klein}): 
(i) 
a regular dihedral group $\mathcal{D}_{2n}$ with $n>1$; 
(ii) 
the alternating group $\mathcal{A}_4$ of degree $4$; 
(iii) 
the symmetric group   $\mathcal{S}_4$ of degree $4$; 
(iv) 
the alternating group $\mathcal{A}_5$ of degree $5$ 
(notice that $G_{\hat{\nabla} , \varepsilon}$ is not isomorphic to 
a cyclic group $\mbox{\boldmath{$Z$}} /n\mbox{\boldmath{$Z$}}$, 
because of $\phi_{\varepsilon} \in (0, \pi /2]$). 

The following holds: 

\begin{proposition}\label{pro:case(i)} 
If $G_{\hat{\nabla} , \varepsilon}$ is isomorphic to $\mathcal{D}_{2n}$ 
with $n>1$, 
then one of the following holds\/$:$ 
\begin{itemize} 
\item[{\rm (a)}]{$\theta_{\varepsilon , 1} 
                 =\theta_{\varepsilon , 2} =\pi$, 
                 $\phi_{\varepsilon} =p\pi /n$,} 
\item[{\rm (b)}]{$\{ 
                  \theta_{\varepsilon , 1} , 
                  \theta_{\varepsilon , 2} \} =\{ \pi , 2p\pi /n\}$, 
                 $\phi_{\varepsilon} =\pi /2$,} 
\end{itemize} 
where $p$ is an integer such that $n$ and $p$ are relatively prime. 
\end{proposition} 

Moreover, the following hold: 

\begin{proposition}[\cite{AK}]\label{pro:case(ii)} 
If $G_{\hat{\nabla} , \varepsilon}$ is isomorphic to $\mathcal{A}_4$, 
then one of the following holds\/$:$ 
\begin{itemize} 
\item[${\rm (a)}$]{$\theta_{\varepsilon , 1} =\pi$, 
$\theta_{\varepsilon , 2} \in \{ 2\pi /3, 4\pi /3\}$, 
$\cos \phi_{\varepsilon} =\sqrt{1/3}$\,$;$} 
\item[${\rm (b)}$]{$\theta_{\varepsilon , 1} \in \{ 2\pi /3, 4\pi /3\}$ and 
$$(\theta_{\varepsilon , 2} , \cos \phi_{\varepsilon} ) 
   \in \{ ( \pi ,   \sqrt{1/3} ), \ 
          (2\pi /3, 1/3), \ 
          (4\pi /3, 1/3)\} .$$} 
\end{itemize} 
\end{proposition} 

\begin{proposition}[\cite{AK}]\label{pro:case(iii)} 
If $G_{\hat{\nabla} , \varepsilon}$ is isomorphic to $\mathcal{S}_4$, 
then one of the following holds\/$:$ 
\begin{itemize} 
\item[{\rm (a)}]{$\theta_{\varepsilon , 1} =\pi$, and 
$(\theta_{\varepsilon , 2} , \phi_{\varepsilon} )$ satisfies one of 
\begin{itemize} 
\item[{\rm (a-1)}]{$\theta_{\varepsilon , 2} \in \{ 2\pi /3, 4\pi /3\}$, 
                   $\cos \phi_{\varepsilon} =\sqrt{2/3}$,} 
\item[{\rm (a-2)}]{$\theta_{\varepsilon , 2} \in \{ \pi /2, 3\pi /2\}$, 
                   $\phi_{\varepsilon} =\pi /4;$} 
\end{itemize}} 
\item[{\rm (b)}]{$\theta_{\varepsilon , 1} \in \{ 2\pi /3, 4\pi /3\}$, and 
$$(\theta_{\varepsilon , 2} , \cos \phi_{\varepsilon} ) 
\in \{ ( \pi   , \sqrt{2/3} ), \ 
       ( \pi /2, \sqrt{1/3} ), \ 
       (3\pi /2, \sqrt{1/3} )\} ;$$} 
\item[{\rm (c)}]{$\theta_{\varepsilon , 1} \in \{ \pi /2, 3\pi /2\}$, and 
$(\theta_{\varepsilon , 2} , \phi_{\varepsilon} )$ satisfies one of 
\begin{itemize} 
\item[{\rm (c-1)}]{$\theta_{\varepsilon , 2} \in \{ 2\pi /3, 4\pi /3\}$, 
                   $\cos \phi_{\varepsilon} =\sqrt{1/3}$,} 
\item[{\rm (c-2)}]{$(\theta_{\varepsilon , 2} , \phi_{\varepsilon} ) 
                     \in \{ ( \pi   , \pi /4), 
                            ( \pi /2, \pi /2), 
                            (3\pi /2, \pi /2)\}$.} 
\end{itemize}} 
\end{itemize} 
\end{proposition} 

We set $\rho :=(1/2)(1+\sqrt{5} )$ and 
\begin{equation*} 
\begin{split} 
& \cos \phi_{(2, 3)}    =\dfrac{\rho}{\sqrt{3}} , \quad 
  \cos \phi_{(2, 5), k} =\sqrt{\dfrac{\rho^{3-2k}}{\sqrt{5}}} \quad (k=1, 2), 
  \\ 
& \cos \phi_{(3, 3)}    =\dfrac{\sqrt{5}}{3} , \quad 
  \cos \phi_{(3, 5), k} =\sqrt{\dfrac{\rho^{3(3-2k)}}{3\sqrt{5}}} 
                         \quad (k=1, 2), \\ 
& \cos \phi_{(5, 5)}    =\dfrac{1}{\sqrt{5}} . 
\end{split} 
\end{equation*} 
The following holds: 

\begin{proposition}[\cite{AK}]\label{pro:case(iv)} 
If $G_{\hat{\nabla} , \varepsilon}$ is isomorphic to $\mathcal{A}_5$, 
then one of the following holds\/$:$ 
\begin{itemize} 
\item[{\rm (a)}]{$\theta_{\varepsilon , 1} =\pi$, and 
$(\theta_{\varepsilon , 2} , \phi_{\varepsilon} )$ satisfies one of 
\begin{itemize} 
\item[{\rm (a-1)}]{$\theta_{\varepsilon , 2} \in \{ 2\pi /3, 4\pi /3\}$, 
                   $\phi_{\varepsilon} =\phi_{(2, 3)}$,} 
\item[{\rm (a-2)}]{$\theta_{\varepsilon , 2} 
                    \in \{ 2\pi /5, 4\pi /5, 6\pi /5, 8\pi /5\}$, 
                   $\phi_{\varepsilon} \in \{ \phi_{(2, 5), 1} , 
                                              \phi_{(2, 5), 2}\} ;$} 
\end{itemize}} 
\item[{\rm (b)}]{$\theta_{\varepsilon , 1} \in \{ 2\pi /3, 4\pi /3\}$, and 
$(\theta_{\varepsilon , 2} , \phi_{\varepsilon} )$ satisfies one of 
\begin{itemize} 
\item[{\rm (b-1)}]{$(\theta_{\varepsilon , 2} , \phi_{\varepsilon} ) 
                     \in \{ ( \pi ,   \phi_{(2, 3)} ), \ 
                            (2\pi /3, \phi_{(3, 3)} ), \ 
                            (4\pi /3, \phi_{(3, 3)} )\}$,} 
\item[{\rm (b-2)}]{$\theta_{\varepsilon , 2} 
                    \in \{ 2\pi /5, 4\pi /5, 6\pi /5, 8\pi /5\}$, 
                   $\phi_{\varepsilon} \in \{ \phi_{(3, 5), 1} , 
                                              \phi_{(3, 5), 2}\} ;$} 
\end{itemize}} 
\item[{\rm (c)}]{$\theta_{\varepsilon , 1} 
                  \in \{ 2\pi /5, 4\pi /5, 6\pi /5, 8\pi /5\}$, 
and $(\theta_{\varepsilon , 2} , \phi_{\varepsilon} )$ satisfies one of 
\begin{itemize} 
\item[{\rm (c-1)}]{$(\theta_{\varepsilon , 2} , \phi_{\varepsilon} ) 
                     \in \{ ( \pi , \phi_{(2, 5), 1} ), \ 
                            ( \pi , \phi_{(2, 5), 2} )\}$,} 
\item[{\rm (c-2)}]{$\theta_{\varepsilon , 2} \in \{ 2\pi /3, 4\pi /3\}$, 
                   $\phi_{\varepsilon} \in \{ \phi_{(3, 5), 1} , 
                                              \phi_{(3, 5), 2} \}$,}
\item[{\rm (c-3)}]{$\theta_{\varepsilon , 2} 
                        \in \{ 2\pi /5, 4\pi /5, 6\pi /5, 8\pi /5\}$, 
$\phi_{\varepsilon} =\phi_{(5, 5)}$.} 
\end{itemize}} 
\end{itemize} 
\end{proposition} 

By the double covering $SO(4)\longrightarrow SO(3)\times SO(3)$, 
we have 

\begin{proposition}\label{pro:finite} 
The topological holonomy group $G_{\nabla}$ of $\nabla$ is finite 
if and only if 
both of the topological holonomy groups $G_{\hat{\nabla} , \pm}$ are 
finite. 
\end{proposition} 

\section{Topological holonomy groups 
with generators of order two}\label{sect:ord2thg} 

\setcounter{equation}{0} 

The following is basic in this paper. 

\begin{proposition}\label{pro:pipi/2} 
Suppose $\{ {\rm ord}\,(C_{\varepsilon , 1} ), 
            {\rm ord}\,(C_{\varepsilon , 2} )\} =\{ 2, \infty \}$ 
and     $\phi_{\varepsilon} =\pi /2$. 
Then $G_{\hat{\nabla} , \varepsilon}$ is infinite 
but not dense in $SO(3)$. 
\end{proposition} 

\begin{remark}\label{rem:A} 
Let $C_{\varepsilon , 1}$, $C_{\varepsilon , 2}$ be 
as in Proposition~\ref{pro:pipi/2}. 
Then there exists an element $\omega_0$ 
of $\hat{E}_{\varepsilon , {\rm pr}(0,0)}$ satisfying $d(\omega_0 )=2$, 
and for an element $\omega$ 
of $\hat{E}_{\varepsilon , {\rm pr}(0,0)} \setminus \{ \pm \omega_0 \}$, 
$X(\omega )$ is given by dense subsets of one or two circles 
in $\hat{E}_{\varepsilon , {\rm pr}(0,0)}$. 
\end{remark} 

We set 
\begin{equation*} 
U(\phi ) 
:=\left[ \begin{array}{ccc} 
         \cos \phi & \sin \phi & 0 \\ 
        -\sin \phi & \cos \phi & 0 \\ 
          0        &  0        & 1 
           \end{array} 
  \right] 
  \quad 
 (\phi \in (0, \pi /2 ]). 
\end{equation*} 
We will prove 

\begin{proposition}\label{pro:pipi} 
Suppose that $C_{\varepsilon , 1}$, $C_{\varepsilon , 2}$ are 
of order two and not commutative. 
If $\phi_{\varepsilon} \in (0, \pi /2)$ is represented 
as $q\pi$ for a rational number $q$, 
then $G_{\hat{\nabla} , \varepsilon}$ is finite\/$;$ 
if $\phi_{\varepsilon}$ is represented as $\psi \pi$ 
for an irrational number $\psi$, 
then $G_{\hat{\nabla} , \varepsilon}$ is conjugate to 
a subgroup of $SO(3)$ given in Proposition~\ref{pro:pipi/2}, 
and therefore infinite but not dense in $SO(3)$. 
\end{proposition} 

\vspace{3mm} 

\par\noindent 
\textit{Proof} \ 
We can set $C_{\varepsilon , 1} :=        C(\pi )$, 
           $C_{\varepsilon , 2} := U(\phi_{\varepsilon} )C(\pi ){}^t 
                                   U(\phi_{\varepsilon} )$. 
Then we have $C_{\varepsilon , 1} C_{\varepsilon , 2} 
                          =U(-2\phi_{\varepsilon} )$, 
             $C_{\varepsilon , 2} C_{\varepsilon , 1} 
                          =U( 2\phi_{\varepsilon} )$. 
Therefore $G_{\hat{\nabla} , \varepsilon}$ is generated 
by  $C_{\varepsilon , 1}$ 
and $C_{\varepsilon , 2} C_{\varepsilon , 1}$. 
The angle between the rotation axes 
of $C_{\varepsilon , 1}$, $C_{\varepsilon , 2} C_{\varepsilon , 1}$ 
is equal to $\pi /2$. 
Therefore, if $\phi_{\varepsilon} \in (0, \pi /2)$ is represented as $q\pi$, 
then $G_{\hat{\nabla} , \varepsilon}$ is finite, 
and if $\phi_{\varepsilon} =\psi \pi$, 
then $G_{\hat{\nabla} , \varepsilon}$ is conjugate to 
a subgroup of $SO(3)$ given in Proposition~\ref{pro:pipi/2}. 
Therefore by Proposition~\ref{pro:pipi/2}, 
we obtain Proposition~\ref{pro:pipi}. 
\hfill 
$\square$ 

\vspace{3mm} 

\begin{remark} 
Let $C_{\varepsilon , 1}$, $C_{\varepsilon , 2}$ be as in 
Proposition~\ref{pro:pipi}. 
Let $\Gamma$ be the great circle of the unit $2$-sphere $S^2$ 
given by the plane $\Pi$ in $\mbox{\boldmath{$R$}}^3$ 
spanned by eigenvectors 
of $C_{\varepsilon , 1}$, $C_{\varepsilon , 2}$ 
corresponding to $\lambda =1$. 
Then for each point $p$ of $\Gamma$, 
$X(p)$ is finite or dense in $\Gamma$ 
according to $\phi_{\varepsilon} =q\pi$ or $\psi \pi$. 
A nonzero vector orthogonal to the plane $\Pi$ is an eigenvector 
of any element of $G_{\hat{\nabla} , \varepsilon}$. 
\end{remark} 

\begin{remark} 
If $C_{\varepsilon , 1}$, $C_{\varepsilon , 2}$ are of order two and 
if $\phi_{\varepsilon} =\pi /2$, 
then these are commutative. 
\end{remark} 

We will prove 

\begin{proposition}\label{pro:2irr} 
Suppose $\{ {\rm ord}\,(C_{\varepsilon , 1} ), 
            {\rm ord}\,(C_{\varepsilon , 2} )\} =\{ 2, \infty \}$ 
and     $\phi_{\varepsilon} \in (0, \pi /2)$. 
then the topological holonomy group $G_{\hat{\nabla} , \varepsilon}$ has 
a noncommutative pair $(\hat{C}_1 , \hat{C}_2 )$ 
with $\infty ={\rm ord}\,(\hat{C}_1 ) 
      \geq    {\rm ord}\,(\hat{C}_2 ) \geq 3$. 
\end{proposition} 

\vspace{3mm} 

\par\noindent 
\textit{Proof} \ 
We can suppose $C_{\varepsilon , 1} 
          =U(\phi_{\varepsilon} )C(\pi ){}^t 
           U(\phi_{\varepsilon} )$ 
for $\phi_{\varepsilon} \in (0, \pi /2)$ 
and $C_{\varepsilon , 2} =C(\psi \pi )$ 
for an irrational number $\psi$. 
Then for $p_1 ={}^t [1 \ 0 \ 0]$, we have 
\begin{equation*} 
\begin{split} 
C_{\varepsilon , 1} C_{\varepsilon , 2} p_1 
& 
=\left[ \begin{array}{c} 
         \cos 2\phi_{\varepsilon} \\ 
        -\sin 2\phi_{\varepsilon} \\ 
                0 
          \end{array} 
 \right] , \\  
C_{\varepsilon , 2} C_{\varepsilon , 1} C_{\varepsilon , 2} p_1 
& 
=\left[ \begin{array}{c} 
         \cos 2\phi_{\varepsilon}               \\ 
        -\sin 2\phi_{\varepsilon} \cos \psi \pi \\ 
        -\sin 2\phi_{\varepsilon} \sin \psi \pi 
          \end{array} 
 \right] 
\end{split} 
\end{equation*} 
and therefore noticing $\phi_{\varepsilon} \in (0, \pi /2)$ 
and $\psi \in \mbox{\boldmath{$R$}} \setminus \mbox{\boldmath{$Q$}}$, 
we see that the third component 
of $(C_{\varepsilon , 1} C_{\varepsilon , 2} )^2 p_1$ is nonzero. 
This means that the order 
of $C_{\varepsilon , 1} C_{\varepsilon , 2}$ is infinite, 
or finite and more than two and 
that eigenvectors 
of $C_{\varepsilon , 1} C_{\varepsilon , 2}$ 
and $C_{\varepsilon , 2}$ corresponding to $\lambda =1$ 
are linearly independent. 
Therefore $\hat{C}_1 :=C_{\varepsilon , 1} C_{\varepsilon , 2}$, 
$\hat{C}_2 :=C_{\varepsilon , 2}$ are noncommutative 
and satisfy $\infty ={\rm ord}\,(\hat{C}_1 ) 
             \geq    {\rm ord}\,(\hat{C}_2 ) \geq 3$. 
Hence we obtain Proposition~\ref{pro:2irr}. 
\hfill 
$\square$ 

\vspace{3mm} 

\begin{remark} 
Notice that in Proposition~\ref{pro:2irr}, 
we can not remove $\phi_{\varepsilon} \in (0, \pi /2)$, that is, 
we can not set $\phi_{\varepsilon} :=\pi /2$ 
(see Proposition~\ref{pro:pipi/2}). 
\end{remark} 

\section{Topological holonomy groups 
with generators of order more than two}\label{sect:ord>2thg} 

\setcounter{equation}{0} 

The following holds: 

\begin{proposition}\label{pro:cond1} 
For $\varepsilon \in \{ +, -\}$, suppose 
\begin{itemize} 
\item[{\rm (a)}]{one of 
$\theta_{\varepsilon , 1} /\pi$, $\theta_{\varepsilon , 2} /\pi$ is 
an irrational number, and} 
\item[{\rm (b)}]{the other is not an integer.} 
\end{itemize} 
Then the topological holonomy group $G_{\hat{\nabla} , \varepsilon}$ has 
a noncommutative pair $(\hat{C}_1 , \hat{C}_2 )$ 
with $\infty ={\rm ord}\,(\hat{C}_1 ) 
      \geq    {\rm ord}\,(\hat{C}_2 ) \geq 3$. 
\end{proposition} 

\vspace{3mm} 

\par\noindent 
\textit{Proof} \ 
If we suppose ${\rm ord}\,(C_{\varepsilon , 1} )=\infty$, 
then $\hat{C}_1 :=C_{\varepsilon , 1}$, 
     $\hat{C}_2 :=C_{\varepsilon , 2}$ satisfy 
     $\infty ={\rm ord}\,(\hat{C}_1 ) 
      \geq    {\rm ord}\,(\hat{C}_2 ) \geq 3$. 
Noticing $\phi_{\varepsilon} \in (0, \pi /2]$, 
we see that $\hat{C}_1$, $\hat{C}_2$ are noncommutative. 
Hence we obtain Proposition~\ref{pro:cond1}. 
\hfill 
$\square$ 

\vspace{3mm} 

One of the generators of topological holonomy groups 
as in Propositions~\ref{pro:2irr} and \ref{pro:cond1} has 
infinite order. 
In the following, 
we will find examples of topological holonomy groups 
in $\bigwedge^2_{\varepsilon}\!E$ 
generated by two finite order elements 
and equipped with noncommutative pairs 
which consist of infinite order elements. 

Suppose $(\theta_{\varepsilon , 1} , 
          \theta_{\varepsilon , 2} , 
          \phi_{\varepsilon} ) 
        =(\pi , 2\pi /n , \pi /2)$ for a prime number $n>2$ or 
that $(\theta_{\varepsilon , 1} , 
       \theta_{\varepsilon , 2} , 
       \phi_{\varepsilon} )$ is a triplet given in one of 
Propositions~\ref{pro:case(ii)}, \ref{pro:case(iii)}, \ref{pro:case(iv)}. 
Let $(\theta'_{\varepsilon , 1} , 
      \theta'_{\varepsilon , 2} )$ be 
either $(\theta_{\varepsilon , 1} /2, 
         \theta_{\varepsilon , 2}   )$ 
or     $(\theta_{\varepsilon , 1}   , 
         \theta_{\varepsilon , 2} /2)$ 
for    $(\theta_{\varepsilon , 1}   , 
         \theta_{\varepsilon , 2}   )$ 
and set 
\begin{equation} 
     C'_{\varepsilon , 1} 
   :=C(\theta'_{\varepsilon , 1} ), \quad 
     C'_{\varepsilon , 2} 
   :=U(\phi_{\varepsilon} )C(\theta'_{\varepsilon , 2} )
{}^t U(\phi_{\varepsilon} ). 
\label{C'}
\end{equation} 
In the following, 
we suppose 
\begin{equation*} 
0<\theta'_{\varepsilon , 1} \leq \pi , \quad 
0<\theta'_{\varepsilon , 2} \leq \pi . 
\end{equation*} 
Let $\nabla$ be an $h$-connection of $E$ satisfying \eqref{CxCy} 
with $C_{\varepsilon , 1} :=C'_{\varepsilon , 1}$, 
     $C_{\varepsilon , 2} :=C'_{\varepsilon , 2}$ 
(notice Proposition~\ref{pro:prod2}). 
Then the topological holonomy group $G_{\hat{\nabla} , \varepsilon}$ 
of $\hat{\nabla}$ 
in $\bigwedge^2_{\varepsilon}\!E$ at ${\rm pr}(0,0)$ is generated by 
finite order elements $C'_{\varepsilon , 1}$, $C'_{\varepsilon , 2}$. 
We will prove 

\begin{proposition}\label{pro:inftyinfty} 
If $G_{\hat{\nabla} , \varepsilon}$ is infinite, 
then two elements $C'_{\varepsilon , 1} C'_{\varepsilon , 2}$, 
                  $C'_{\varepsilon , 2} C'_{\varepsilon , 1}$ 
are noncommutative and of infinite order. 
\end{proposition} 

In order to prove Proposition~\ref{pro:inftyinfty}, 
we need lemmas. 
The following holds: 

\begin{lemma}\label{lem:evls} 
The set of eigenvalues of $C'_{\varepsilon , 1} 
                           C'_{\varepsilon , 2}$ 
coincides with 
the set of eigenvalues of $C'_{\varepsilon , 2} 
                           C'_{\varepsilon , 1}$. 
\end{lemma} 

We will prove 

\begin{lemma}\label{lem:1evcts} 
An eigenvector of $C'_{\varepsilon , 1} 
                   C'_{\varepsilon , 2}$ 
corresponding to $\lambda =1$ and 
an eigenvector of $C'_{\varepsilon , 2} 
                   C'_{\varepsilon , 1}$ 
corresponding to $\lambda =1$ are linearly independent. 
\end{lemma} 

\vspace{3mm} 

\par\noindent 
\textit{Proof} \ 
Let $\eta$, $\eta^*$ be eigenvectors 
of $C'_{\varepsilon , 1} C'_{\varepsilon , 2}$, 
   $C'_{\varepsilon , 2} C'_{\varepsilon , 1}$ 
respectively corresponding to $\lambda =1$. 
We can choose $\eta$, $\eta^*$ so that they are unit vectors. 
Suppose that  $\eta$, $\eta^*$ are linearly dependent. 
Then we have $\eta^* =\eta$ or $-\eta$. 
Since   $(I_3 -C'_{\varepsilon , 1} C'_{\varepsilon , 2} )\eta =0$, 
we have $(I_3 -C'_{\varepsilon , 2} C'_{\varepsilon , 1} ) 
                               {}^t C'_{\varepsilon , 1} \eta =0$. 
This means $\eta^* =\pm {}^t C'_{\varepsilon , 1} \eta$. 
Therefore we obtain $C'_{\varepsilon , 1} \eta =\eta$ or $-\eta$. 
If $C'_{\varepsilon , 1} \eta =\eta$, 
then $C'_{\varepsilon , 2} C'_{\varepsilon , 1} \eta^* =\eta^*$ 
means $C'_y \eta =\eta$, 
which contradicts $\phi_{\varepsilon} \not= 0$. 
If $C'_{\varepsilon , 1} \eta =-\eta$, 
then $\theta'_{\varepsilon , 1} =\pi$, 
and  $C'_{\varepsilon , 2} \eta =-\eta$, 
which means $\theta'_{\varepsilon , 2} =\pi$. 
However, in this section, we do not study 
the case of $(\theta'_{\varepsilon , 1} , 
              \theta'_{\varepsilon , 2} ) 
            =(\pi , \pi )$. 
Hence we have concluded 
that $\eta$, $\eta^*$ are linearly independent. 
\hfill 
$\square$ 

\vspace{3mm} 

Noticing that $C'_{\varepsilon , 1} C'_{\varepsilon , 2}$ is 
an element of $SO(3)$, 
we see that the characteristic polynomial 
of $C'_{\varepsilon , 1} C'_{\varepsilon , 2}$ is given by 
\begin{equation} 
\begin{split} 
   \chi_{C'_{\varepsilon , 1} C'_{\varepsilon , 2}} (\lambda ) 
& =\det (\lambda E_3 -C'_{\varepsilon , 1} C'_{\varepsilon , 2} ) \\ 
& =  \lambda^3 
  -({\rm tr}\,(C'_{\varepsilon , 1} C'_{\varepsilon , 2} )) 
     \lambda^2 
  +({\rm tr}\,(C'_{\varepsilon , 1} C'_{\varepsilon , 2} )) 
     \lambda -1 \\ 
& =(\lambda   -1) 
   (\lambda^2 
  +(1-({\rm tr}\,(C'_{\varepsilon , 1} C'_{\varepsilon , 2} ))) 
    \lambda   
  + 1). 
\end{split} 
\label{cp}
\end{equation} 
By \eqref{C'}, we obtain 
\begin{equation} 
\begin{split} 
&    {\rm tr}\,(C'_{\varepsilon , 1} C'_{\varepsilon , 2} ) \\ 
& =   \cos   \theta'_{\varepsilon , 1}  
      \cos   \theta'_{\varepsilon , 2} 
  - 2 \sin   \theta'_{\varepsilon , 1}  
      \sin   \theta'_{\varepsilon , 2}  
      \cos   \phi_{\varepsilon} \\ 
&     \quad 
  +(1+\cos   \theta'_{\varepsilon , 1}  
      \cos   \theta'_{\varepsilon , 2} ) 
      \cos^2 \phi_{\varepsilon} 
  +(  \cos   \theta'_{\varepsilon , 1} 
     +\cos   \theta'_{\varepsilon , 2} ) 
      \sin^2 \phi_{\varepsilon} . 
\end{split}
\label{trace}
\end{equation} 

Let $\zeta$ be an eigenvalue 
of $C'_{\varepsilon , 1} C'_{\varepsilon , 2}$ other than $1$. 
We will prove 

\begin{lemma}\label{lem:case(ii)infty}
Let $(\theta_{\varepsilon , 1} , 
      \theta_{\varepsilon , 2} , 
      \phi_{\varepsilon} )$ be 
one of the triplets given in Proposition~\ref{pro:case(ii)} 
and suppose that $G_{\hat{\nabla} , \varepsilon}$ is infinite. 
Then $(\theta'_{\varepsilon , 1} , 
       \theta'_{\varepsilon , 2} , 
       \phi_{\varepsilon} )$ satisfies 
$\{ \theta'_{\varepsilon , 1} , 
    \theta'_{\varepsilon , 2} \} =\{ \pi , \pi /3\}$ and 
$\cos \phi_{\varepsilon} =\sqrt{1/3}$, 
and $\zeta$ is represented as $\zeta =\exp (\sqrt{-1} \psi \pi )$ 
for a real, irrational number $\psi$. 
\end{lemma} 

Related to the proof of Lemma~\ref{lem:case(ii)infty}, 
we recall cyclotomic polynomials and their properties. 
For a natural number $n\in \mbox{\boldmath{$N$}}$, 
let $P(n)$ be the set of primitive $n$-th roots of unity. 
Then the number of elements of $P(n)$ is just 
the value $\varphi (n)$ of Euler's totient function for $n$. 
Let $k_n$ be a polynomial of one variable $\lambda$ defined by 
\begin{equation*} 
k_n (\lambda )=\prod_{\eta \in P(n)} (\lambda -\eta ).  
\end{equation*} 
Then $k_n$ is called the $n$-th \textit{cyclotomic polynomial}. 
Since $\sharp P(n)=\varphi (n)$, 
$k_n$ is a polynomial of degree $\varphi (n)$. 
For a natural number $N\in \mbox{\boldmath{$N$}}$, we have 
\begin{equation*} 
\lambda^N -1 =\prod_{n\in D(N)} k_n (\lambda ), 
\end{equation*} 
where $D(N)$ denotes the set of divisors of $N$. 
By induction, we can show that the coefficients of $k_n$ are integers. 
In addition, $k_n$ is irreducible over $\mbox{\boldmath{$Q$}}$ 
(\cite[Theorem 3.8.6]{nagata}). 

We will prove Lemma~\ref{lem:case(ii)infty}. 
By Proposition~\ref{pro:case(ii)} and direct computations, 
we observe that 
if $(\theta_{\varepsilon , 1} , 
     \theta_{\varepsilon , 2} , 
     \phi_{\varepsilon} )$ is as in Proposition~\ref{pro:case(ii)}, 
then $\{ \theta'_{\varepsilon , 1} , 
         \theta'_{\varepsilon , 2} \} =\{ \pi , \pi /3\}$ and 
$\cos \phi_{\varepsilon} =\sqrt{1/3}$, and 
therefore by \eqref{trace}, 
we have ${\rm tr}\,(C'_{\varepsilon , 1} 
                    C'_{\varepsilon , 2} )=-2/3$. 
In addition, by \eqref{cp}, 
$\zeta$ is algebraic on $\mbox{\boldmath{$Q$}}$ 
and the minimal polynomial $f_{\zeta}$ of $\zeta$ 
on $\mbox{\boldmath{$Q$}}$ is represented as 
\begin{equation*} 
f_{\zeta} (\lambda )=\lambda^2 +\dfrac{5}{3} \lambda +1  
\end{equation*} 
(see \cite[Section 3.1]{nagata} 
for algebraic numbers and their minimal polynomials). 
We have $f_{\zeta} (\zeta )=0$. 
If $f=f(\lambda )$ is a polynomial 
with coefficients in $\mbox{\boldmath{$Q$}}$ 
satisfying $f(\zeta )=0$, 
then there exists a polynomial $g(\lambda )$ 
with coefficients in $\mbox{\boldmath{$Q$}}$ 
satisfying $f(\lambda)=g(\lambda )f_{\zeta} (\lambda )$. 
We represent $\zeta$ as $\zeta =\exp (\sqrt{-1} \psi \pi )$ 
for a real number $\psi$. 
Suppose that $\psi$ is a rational number. 
Then there exists a positive number $p$ 
such that $f(\lambda ):=\lambda^{2p} -1$ satisfies $f(\zeta )=0$. 
Therefore we have $f(\lambda)=g(\lambda )f_{\zeta} (\lambda )$ 
and then $p>1$. 
However, 
since the cyclotomic polynomials are irreducible 
over $\mbox{\boldmath{$Q$}}$ and 
since their coefficients are integers, 
$f$ never has such a representation. 
Therefore $\psi$ is an irrational number and 
we obtain Lemma~\ref{lem:case(ii)infty}. 
On the other hand, 
we can obtain the same conclusion by direct computations as follows. 
We represent $g$ as $g(\lambda )=\sum^{2p-2}_{k=0} A_k \lambda^k$, 
where $A_0 =-1$, $A_k \in \mbox{\boldmath{$Q$}}$. 
Then we have $A_k = B_k /3^k$, 
where $B_k$ is an integer and not devided by $3$. 
This contradicts $A_{2p-2-k} =-A_k$. 
Therefore we see that $\psi$ is an irrational number again. 

We will prove 

\begin{lemma}\label{lem:case(iii)infty}
Let $(\theta_{\varepsilon , 1} , 
      \theta_{\varepsilon , 2} , 
      \phi_{\varepsilon} )$ be 
one of the triplets given in Proposition~\ref{pro:case(iii)} 
and suppose that $G_{\hat{\nabla} , \varepsilon}$ is infinite. 
Then $(\theta'_{\varepsilon , 1} , 
       \theta'_{\varepsilon , 2} , 
       \phi_{\varepsilon} )$ satisfies one of the following\/$:$ 
\begin{itemize} 
\item[{\rm (a)}]{$\phi_{\varepsilon} =\pi /2$, 
and $\{ \theta'_{\varepsilon , 1} , 
        \theta'_{\varepsilon , 2} \} 
   = \{ \pi /2,  \pi /4\}$ 
 or $\{ \pi /2, 3\pi /4\};$} 
\item[{\rm (b)}]{$\phi_{\varepsilon} =\pi /4$, 
and $\{ \theta'_{\varepsilon , 1} , 
        \theta'_{\varepsilon , 2} \}$ is equal to  
    $\{ \pi /2,  \pi /2\}$, 
    $\{ \pi   ,  \pi /4\}$  
 or $\{ \pi   , 3\pi /4\};$} 
\item[{\rm (c)}]{$\cos \phi_{\varepsilon} =\sqrt{1/3}$,  
and $\{ \theta'_{\varepsilon , 1} , 
        \theta'_{\varepsilon , 2} \}$ is equal to  
    $\{  \pi /2,  \pi /3\}$, 
    $\{  \pi /4, 2\pi /3\}$  
 or $\{ 3\pi /4, 2\pi /3\};$} 
\item[{\rm (d)}]{$\cos \phi_{\varepsilon} =\sqrt{2/3}$, 
and $\{ \theta'_{\varepsilon , 1} , 
        \theta'_{\varepsilon , 2} \} 
   = \{ \pi   ,  \pi /3\}$  
 or $\{ \pi /2, 2\pi /3\}$.} 
\end{itemize} 
In any case, 
$\zeta$ is represented as $\zeta =\exp (\sqrt{-1} \psi \pi )$ 
for a real, irrational number $\psi$. 
\end{lemma} 

\vspace{3mm} 

\par\noindent 
\textit{Proof} \ 
By Proposition~\ref{pro:case(iii)} and direct computations, 
if $(\theta_{\varepsilon , 1} , 
     \theta_{\varepsilon , 2} , 
     \phi_{\varepsilon} )$ is as in Proposition~\ref{pro:case(iii)}, 
then $\{ \theta'_{\varepsilon , 1} , 
         \theta'_{\varepsilon , 2} \}$ and $\phi_{\varepsilon}$ 
satisfy (a), (b), (c) or (d) in Lemma~\ref{lem:case(iii)infty}. 
If $(\{ \theta'_{\varepsilon , 1} , 
        \theta'_{\varepsilon , 2} \} , 
           \phi_{\varepsilon} )$ coincides with one of 
\begin{equation*} 
(\{ \pi   ,  \pi /3\} , \cos^{-1} (\sqrt{2/3} )), \ 
(\{ \pi /2,  \pi /3\} , \cos^{-1} (\sqrt{1/3} )), 
\end{equation*} 
then by \eqref{trace}, 
we have ${\rm tr}\,(C'_{\varepsilon , 1} 
                    C'_{\varepsilon , 2} )=-1/3$. 
In addition, by \eqref{cp}, 
$\zeta$ is algebraic on $\mbox{\boldmath{$Q$}}$ 
and $f_{\zeta}$ is represented as 
\begin{equation*} 
f_{\zeta} (\lambda )=\lambda^2 +\dfrac{4}{3} \lambda +1 . 
\end{equation*} 
Therefore referring to the proof of Lemma~\ref{lem:case(ii)infty}, 
we observe 
that $\zeta$ is represented as $\zeta =\exp (\sqrt{-1} \psi \pi )$ 
for a real, irrational number $\psi$. 
If 
\begin{equation*} 
(\{ \theta'_{\varepsilon , 1} , 
    \theta'_{\varepsilon , 2} \} , 
    \phi_{\varepsilon} )
 \in 
\{ (\{ \pi /2, 2\pi /3\} , \cos^{-1} (\sqrt{2/3} )), 
   (\{ \pi /2,  \pi /2\} , \pi /4 ) 
\} , 
\end{equation*} 
then by \eqref{trace}, 
we have ${\rm tr}\,(C'_{\varepsilon , 1} 
                    C'_{\varepsilon , 2} )=1/2 -\sqrt{2}$, 
and by \eqref{cp}, 
$\zeta$ is algebraic on $\mbox{\boldmath{$Q$}}$ 
and $f_{\zeta}$ is represented as 
\begin{equation} 
f_{\zeta} (\lambda )
=\lambda^4 +\lambda^3 +\dfrac{1}{4} \lambda^2 +\lambda +1; 
\label{deg4_1} 
\end{equation} 
if $(\{ \theta'_{\varepsilon , 1} , 
        \theta'_{\varepsilon , 2} \} , 
           \phi_{\varepsilon} )$ is equal to one of 
\begin{equation*} 
\begin{split} 
&    (\{  \pi   ,  \pi /4\} , \pi /4 ), \ 
     (\{  \pi   , 3\pi /4\} , \pi /4 ), \\ 
&    (\{ 2\pi /3,  \pi /4\} , \cos^{-1} (\sqrt{1/3} )), \ 
     (\{ 2\pi /3, 3\pi /4\} , \cos^{-1} (\sqrt{1/3} )), \\ 
&    (\{  \pi /2,  \pi /4\} , \pi /2 ), \ 
     (\{  \pi /2, 3\pi /4\} , \pi /2 ), 
\end{split} 
\end{equation*} 
then by \eqref{trace}, 
we have ${\rm tr}\,(C'_{\varepsilon , 1} 
                    C'_{\varepsilon , 2} )=1/\sqrt{2}$ 
                                      or $-1/\sqrt{2}$, 
and by \eqref{cp}, 
$\zeta$ is algebraic on $\mbox{\boldmath{$Q$}}$ 
and $f_{\zeta}$ is represented as 
\begin{equation} 
f_{\zeta} (\lambda )
=\lambda^4 +2\lambda^3 +\dfrac{5}{2} \lambda^2 +2\lambda +1. 
\label{deg4_2} 
\end{equation} 
If $f_{\zeta}$ is represented 
as in \eqref{deg4_1} or \eqref{deg4_2}, 
then as in the above proof of Lemma~\ref{lem:case(ii)infty}, 
based on the properties of cyclotomic polynomials, 
we observe 
that $\zeta$ is represented as $\zeta =\exp (\sqrt{-1} \psi \pi )$ 
for a real, irrational number $\psi$. 
On the other hand, 
we can obtain the same conclusion by direct computations as follows. 
Suppose that $\psi$ is a rational number. 
Then there exists a positive integer $p$ 
such that $f(\lambda ):=\lambda^{2p} -1$ satisfies $f(\zeta )=0$. 
Therefore we have $f(\lambda)=g(\lambda )f_{\zeta} (\lambda )$ 
and $p>2$. 
We represent $g$ as $g(\lambda )=\sum^{2p-4}_{k=0} A_k \lambda^k$, 
where $A_0 =-1$, $A_k \in \mbox{\boldmath{$Q$}}$. 
Suppose \eqref{deg4_2}. 
Then we have $A_{2k} = B_{2k} /2^k$, where $B_{2k}$ is 
an odd integer. 
This contradicts $A_{2p-4 -k} =-A_k$. 
Even if we suppose \eqref{deg4_1}, 
we have a similar contradiction. 
Therefore we see that $\psi$ is an irrational number again. 
Hence we obtain Lemma~\ref{lem:case(iii)infty}. 
\hfill 
$\square$ 

\vspace{3mm} 

We will prove 

\begin{lemma}\label{lem:case(iv(a))infty}
Let $(\theta_{\varepsilon , 1} , 
      \theta_{\varepsilon , 2} , 
      \phi_{\varepsilon} )$ be 
one of the triplets given in Proposition~\ref{pro:case(iv)} 
and suppose that 
one of $\theta_{\varepsilon , 1}$, $\theta_{\varepsilon , 2}$ is 
equal to $\pi$ 
and that $G_{\hat{\nabla} , \varepsilon}$ is infinite. 
Then $(\theta'_{\varepsilon , 1} , 
       \theta'_{\varepsilon , 2} , 
       \phi_{\varepsilon} )$ satisfies one of the following\/$:$ 
\begin{itemize} 
\item[{\rm (a)}]{$\phi_{\varepsilon} =\phi_{(2, 3)}$ 
and $\{ \theta'_{\varepsilon , 1} , 
        \theta'_{\varepsilon , 2} \} 
     \in 
     \{ \{ \pi /2, 2\pi /3\} , 
        \{ \pi   ,  \pi /3\} 
     \} ;$} 
\item[{\rm (b)}]{$\phi_{\varepsilon} =\phi_{(2, 5), 1}$ 
or $\phi_{(2, 5), 2}$, 
and 
\begin{equation*} 
\{ \theta'_{\varepsilon , 1} , 
   \theta'_{\varepsilon , 2} \} 
\in 
\{ \{ \pi   ,  \pi /5\} , 
   \{ \pi   , 3\pi /5\} , 
   \{ \pi /2, 2\pi /5\} , 
   \{ \pi /2, 4\pi /5\} 
\} .
\end{equation*}} 
\end{itemize} 
In any case, 
$\zeta$ is represented as $\zeta =\exp (\sqrt{-1} \psi \pi )$ 
for a real, irrational number $\psi$. 
\end{lemma} 

\vspace{3mm} 

\par\noindent 
\textit{Proof} \ 
By Proposition~\ref{pro:case(iv)} and direct computations, 
if $(\theta_{\varepsilon , 1} , 
     \theta_{\varepsilon , 2} , 
     \phi_{\varepsilon} )$ is as in Lemma~\ref{lem:case(iv(a))infty},  
then $\{ \theta'_{\varepsilon , 1} , 
         \theta'_{\varepsilon , 2} \}$ and $\phi_{\varepsilon}$ 
satisfy (a) or (b) in Lemma~\ref{lem:case(iv(a))infty}. 
If 
\begin{equation*} 
(\{ \theta'_{\varepsilon , 1} , 
    \theta'_{\varepsilon , 2} \} , 
    \phi_{\varepsilon} )
 \in 
\{ (\{ \pi /2, 4\pi /5\} , \phi_{(2, 5), 1} ), 
   (\{ \pi /2, 2\pi /5\} , \phi_{(2, 5), 2} )\} , 
\end{equation*} 
then by \eqref{trace}, 
we have ${\rm tr}\,(C'_{\varepsilon , 1} 
                    C'_{\varepsilon , 2} )=-1/2$. 
In addition, by \eqref{cp}, 
$\zeta$ is algebraic on $\mbox{\boldmath{$Q$}}$ 
and $f_{\zeta}$ is represented as 
\begin{equation*} 
f_{\zeta} (\lambda )=\lambda^2 +\dfrac{3}{2} \lambda +1 . 
\end{equation*} 
Therefore referring to the proof of Lemma~\ref{lem:case(ii)infty}, 
we observe 
that $\zeta$ is represented as $\zeta =\exp (\sqrt{-1} \psi \pi )$ 
for a real, irrational number $\psi$. 
If 
\begin{equation*} 
(\{ \theta'_{\varepsilon , 1} , 
    \theta'_{\varepsilon , 2} \} , 
    \phi_{\varepsilon} )
 \in 
\{ (\{ \pi /2, 2\pi /3\} , \phi_{(2, 3)} ), 
   (\{ \pi /2, 2\pi /5\} , \phi_{(2, 5), 1} ) 
\} , 
\end{equation*} 
then by \eqref{trace}, 
we have ${\rm tr}\,(C'_{\varepsilon , 1} 
                    C'_{\varepsilon , 2} )=(-1/2)\rho$, 
and by \eqref{cp}, 
$\zeta$ is algebraic on $\mbox{\boldmath{$Q$}}$ 
and $f_{\zeta}$ is represented as 
\begin{equation} 
f_{\zeta} (\lambda )
=\lambda^4 +\dfrac{5}{2} \lambda^3 +\dfrac{13}{4} \lambda^2 
           +\dfrac{5}{2} \lambda +1; 
\label{deg4_1'} 
\end{equation} 
if $\{ \theta'_{\varepsilon , 1} , 
       \theta'_{\varepsilon , 2} \} 
   =\{ \pi   ,  \pi /3\}$ and 
   $\phi_{\varepsilon} =\phi_{(2, 3)}$, 
then by \eqref{trace}, 
${\rm tr}\,(C'_{\varepsilon , 1} 
            C'_{\varepsilon , 2} )=(-1/3)\rho^{-2}$, 
and by \eqref{cp}, 
$\zeta$ is algebraic on $\mbox{\boldmath{$Q$}}$ 
and $f_{\zeta}$ is represented as 
\begin{equation} 
f_{\zeta} (\lambda )
=\lambda^4 +3\lambda^3 +\dfrac{37}{9} \lambda^2 +3\lambda +1; 
\label{deg4_2'} 
\end{equation} 
if $(\{ \theta'_{\varepsilon , 1} , 
        \theta'_{\varepsilon , 2} \} , 
        \phi_{\varepsilon} ) 
   =(\{ \pi ,  \pi /5\} , \phi_{(2, 5), 1} )$ 
or $(\{ \pi , 3\pi /5\} , \phi_{(2, 5), 2} )$, 
then by \eqref{trace}, we have 
${\rm tr}\,(C'_{\varepsilon , 1} 
            C'_{\varepsilon , 2} )=-(1/\sqrt{5} )\rho$ 
or                                $-(1/\sqrt{5} )\rho^{-1}$, 
and by \eqref{cp}, 
$\zeta$ is algebraic on $\mbox{\boldmath{$Q$}}$ 
and $f_{\zeta}$ is represented as 
\begin{equation} 
f_{\zeta} (\lambda )
=\lambda^4 +3\lambda^3 +\dfrac{21}{5} \lambda^2 +3\lambda +1; 
\label{deg4_3'} 
\end{equation} 
if $(\{ \theta'_{\varepsilon , 1} , 
        \theta'_{\varepsilon , 2} \} , 
        \phi_{\varepsilon} ) 
   =(\{ \pi , 3\pi /5\} , \phi_{(2, 5), 1} )$ 
or $(\{ \pi ,  \pi /5\} , \phi_{(2, 5), 2} )$, 
then by \eqref{trace}, 
we have ${\rm tr}\,(C'_{\varepsilon , 1} 
                    C'_{\varepsilon , 2} )=2/\sqrt{5}$ or $-2/\sqrt{5}$, 
and by \eqref{cp}, 
$\zeta$ is algebraic on $\mbox{\boldmath{$Q$}}$ 
and $f_{\zeta}$ is represented as 
\begin{equation} 
f_{\zeta} (\lambda )
=\lambda^4 +2\lambda^3 +\dfrac{11}{5} \lambda^2 +2\lambda +1; 
\label{deg4_4'} 
\end{equation} 
if $\{ \theta'_{\varepsilon , 1} , 
       \theta'_{\varepsilon , 2} \} 
   =\{ \pi , 4\pi /5\}$ and 
   $\phi_{\varepsilon} =\phi_{(2, 5), 2}$, 
then by \eqref{trace}, 
${\rm tr}\,(C'_{\varepsilon , 1} 
            C'_{\varepsilon , 2} )=-(3/2)\rho^{-1}$, 
and by \eqref{cp}, 
$\zeta$ is algebraic on $\mbox{\boldmath{$Q$}}$ 
and $f_{\zeta}$ is represented as 
\begin{equation} 
f_{\zeta} (\lambda )
=\lambda^4 +\dfrac{1}{2} \lambda^3 -\dfrac{3}{4} \lambda^2 
           +\dfrac{1}{2} \lambda +1. 
\label{deg4_5'} 
\end{equation} 
If $f_{\zeta}$ is represented 
as in \eqref{deg4_2'}, \eqref{deg4_3'} or \eqref{deg4_4'}, 
then referring to the proof of Lemma~\ref{lem:case(iii)infty}, 
we observe 
that $\zeta$ is represented as $\zeta =\exp (\sqrt{-1} \psi \pi )$ 
for a real, irrational number $\psi$. 
If $f_{\zeta}$ is represented as in \eqref{deg4_1'} or \eqref{deg4_5'}, 
then as in the above proof of Lemma~\ref{lem:case(ii)infty}, 
based on the properties of cyclotomic polynomials, 
we obtain the same result. 
In addition, we can also do it by direct computations as follows. 
Suppose that $f_{\zeta}$ is represented as in \eqref{deg4_1'}. 
We represent $A_k$ as $A_k =B_k /2^k$ 
for $B_k \in \mbox{\boldmath{$Z$}}$. 
Then $B_0 =-1$, $B_1$ is odd and $B_2$ is even. 
Therefore $B_3$ and $B_4$ are odd. 
Inductively, we see that $B_{3k-1}$ is even, and 
that $B_{3k}$, $B_{3k+1}$ are odd for $k\geq 1$. 
This contradicts $A_{2p-4 -k} =-A_k$. 
Even if we suppose that $f_{\zeta}$ is represented as in \eqref{deg4_5'}, 
we have a similar contradiction. 
Therefore we see that $\psi$ is an irrational number again. 
Hence we obtain Lemma~\ref{lem:case(iv(a))infty}. 
\hfill 
$\square$ 

\vspace{3mm} 

We will prove 

\begin{lemma}\label{lem:case(iv(b))infty}
Let $(\theta_{\varepsilon , 1} , 
      \theta_{\varepsilon , 2} , 
      \phi_{\varepsilon} )$ be 
one of the triplets given in Proposition~\ref{pro:case(iv)} 
and suppose that 
one of $\theta_{\varepsilon , 1}$, $\theta_{\varepsilon , 2}$ is 
equal to $2\pi /3$ or $4\pi /3$ 
and that $G_{\hat{\nabla} , \varepsilon}$ is infinite. 
Then $\zeta$ is represented as $\zeta =\exp (\sqrt{-1} \psi \pi )$ 
for a real, irrational number $\psi$. 
\end{lemma} 

\vspace{3mm} 

\par\noindent 
\textit{Proof} \ 
In this proof, 
we consider only the cases of new values 
for the trace ${\rm tr}\,(C'_{\varepsilon , 1} C'_{\varepsilon , 2} )$. 
Therefore by \eqref{trace}, 
we obtain $\theta'_{\varepsilon , 1}  =\pi /3$, and 
$(\theta'_{\varepsilon , 2} , \phi_{\varepsilon} )$ coincides with one of  
\begin{equation*} 
(2\pi /3                    , \phi_{(3, 3)}    ), 
(4\pi /5                    , \phi_{(3, 5), 1} ), 
(2\pi /5                    , \phi_{(3, 5), 2} ), 
(4\pi /5                    , \phi_{(3, 5), 2} ). 
\end{equation*} 
If  $( \theta'_{\varepsilon , 2} , \phi_{\varepsilon} ) 
    =(2\pi /5                    , \phi_{(3, 5), 2} )$, 
then ${\rm tr}\,(C'_{\varepsilon , 1} C'_{\varepsilon , 2} )=2/3$, 
and $\zeta$ is algebraic on $\mbox{\boldmath{$Q$}}$ 
so that $f_{\zeta}$ is represented as 
\begin{equation*} 
f_{\zeta} (\lambda )=\lambda^2 +\dfrac{1}{3} \lambda +1.  
\end{equation*} 
Therefore referring to the proof of Lemma~\ref{lem:case(ii)infty}, 
we observe 
that $\zeta$ is represented as $\zeta =\exp (\sqrt{-1} \psi \pi )$ 
for a real, irrational number $\psi$. 
If $( \theta'_{\varepsilon , 2} , \phi_{\varepsilon} ) 
   =(2\pi /3                    , \phi_{(3, 3)}    )$ 
or $(4\pi /5                    , \phi_{(3, 5), 1} )$, 
then ${\rm tr}\,(C'_{\varepsilon , 1} C'_{\varepsilon , 2} ) 
   =1-(\sqrt{5} /3) \rho^2$, 
and $\zeta$ is algebraic on $\mbox{\boldmath{$Q$}}$ 
so that $f_{\zeta}$ is represented as 
\begin{equation*} 
f_{\zeta} (\lambda )  
  =\lambda^4 +\dfrac{5}{3} \lambda^3 +\dfrac{13}{9} \lambda^2 
             +\dfrac{5}{3} \lambda   +1. 
\end{equation*} 
Therefore referring to the proof of Lemma~\ref{lem:case(iv(a))infty}, 
we observe 
that $\zeta$ is represented as $\zeta =\exp (\sqrt{-1} \psi \pi )$ 
for a real, irrational number $\psi$. 
If  $( \theta'_{\varepsilon , 2} , \phi_{\varepsilon} ) 
    =(4\pi /5                    , \phi_{(3, 5), 2} )$, 
then ${\rm tr}\,(C'_{\varepsilon , 1} C'_{\varepsilon , 2} )=-(1/3)\rho^2$, 
and $\zeta$ is algebraic on $\mbox{\boldmath{$Q$}}$ 
so that $f_{\zeta}$ is represented as in \eqref{deg4_2'}. 
Hence we obtain Lemma~\ref{lem:case(iv(b))infty}. 
\hfill 
$\square$ 

\vspace{3mm} 

We will prove 

\begin{lemma}\label{lem:case(iv(c))infty}
Let $(\theta_{\varepsilon , 1} , 
      \theta_{\varepsilon , 2} , 
      \phi_{\varepsilon} )$ be 
one of the triplets given in Proposition~\ref{pro:case(iv)} 
and suppose that 
one of $\theta_{\varepsilon , 1}$, $\theta_{\varepsilon , 2}$ is 
equal to $2\pi /5$, $4\pi /5$, $6\pi /5$ or $8\pi /5$ 
and that $G_{\hat{\nabla} , \varepsilon}$ is infinite. 
Then $\zeta$ is represented as $\zeta =\exp (\sqrt{-1} \psi \pi )$ 
for a real, irrational number $\psi$. 
\end{lemma} 

\vspace{3mm} 

\par\noindent 
\textit{Proof} \ 
By \eqref{trace}, 
it is seen that there exist no new values 
for the trace ${\rm tr}\,(C'_x C'_y )$. 
Hence we obtain Lemma~\ref{lem:case(iv(c))infty}. 
\hfill 
$\square$ 

\vspace{3mm} 

We will prove 

\begin{lemma}\label{lem:prime} 
Suppose $(\theta_{\varepsilon , 1} , 
          \theta_{\varepsilon , 2} , 
          \phi_{\varepsilon} ) 
        =(\pi , 2\pi /n , \pi /2)$ 
for a prime number $n>2$ and 
that $G_{\hat{\nabla} , \varepsilon}$ is infinite. 
Then $\zeta =\exp (\sqrt{-1} \psi \pi )$, 
where $\psi$ is a real, irrational number. 
\end{lemma} 

\vspace{3mm} 

\par\noindent 
\textit{Proof} \ 
Let $n>2$ be a prime number. 
In this situation, 
we suppose $(\theta'_{\varepsilon , 1} , 
             \theta'_{\varepsilon , 2} , 
             \phi_{\varepsilon} )
           =(\pi /2, 2\pi /n, \pi /2)$. 
Then ${\rm tr}\,(C'_{\varepsilon , 1} C'_{\varepsilon , 2} ) 
      =\cos (2\pi /n)$ by \eqref{trace}, and 
therefore by \eqref{cp}, we have 
\begin{equation*} 
  \chi_{C'_{\varepsilon , 1} C'_{\varepsilon , 2}} (\lambda ) 
=(\lambda -1) 
  \left( \lambda^2 +\left( 1-\cos \dfrac{2\pi}{n} \right) \lambda +1\right) . 
\end{equation*} 
Therefore $\zeta$ satisfies $\zeta^2 +(1-\cos (2\pi /n))\zeta +1=0$. 
This means $\cos (2\pi /n)=\zeta +1 +1/\zeta$. 
Since $n>2$ is prime, 
$\exp (2\sqrt{-1} \pi /n)$ is algebraic on $\mbox{\boldmath{$Q$}}$ 
and its minimal polynomial has degree $n-1$. 
Then we have $2\sum^{\delta}_{j=1} \cos 2j\pi /n =-1$. 
This means that $\cos (2\pi /n)$ is algebraic on $\mbox{\boldmath{$Q$}}$ 
so that its minimal polynomial has degree $\delta =(n-1)/2$. 
Since $\cos (2\pi /n)=\zeta +1 +1/\zeta$, 
$\zeta$ is algebraic on $\mbox{\boldmath{$Q$}}$ 
so that its minimal polynomial has degree $2\delta =n-1$. 
In general, for a positive integer $m$, 
$\cos mx$ is represented as $\cos mx =T_m (\cos x)$, 
where $T_m (X)$ is a polynomial of one variable $X$ of degree $m$, 
and known as the Chebyshev polynomial of degree $m$. 
If we represent $T_m (X)$ as $T_m (X)=\sum^m_{i=0} a_i X^i$, 
then each $a_i$ is an integer, and in particular, 
noticing $\cos mx +\cos (m-2)x =2\cos x\cos (m-1)x$, 
we obtain $a_m =2^{m-1}$.  
Therefore the minimal polynomial of $\cos (2\pi /n)$ is 
represented as 
\begin{equation*} 
F_{\delta} (\lambda )=\sum^{\delta}_{k=0} b_k \lambda^k , 
\end{equation*} 
where 
\begin{equation*} 
b_{\delta} =1, \ 
b_k =\dfrac{b'_k}{2^{\delta -k}} \ (1\leq k\leq \delta -1), \ 
b_0 =\dfrac{2b'_0 +1}{2^{\delta}} , \  
b'_0 , b'_k \in \mbox{\boldmath{$Z$}} . 
\end{equation*} 
Applying $\cos (2\pi /n)=\zeta +1 +1/\zeta$ 
to $F_{\delta} (\cos (2\pi /n))=0$, 
we see that the minimal polynomial of $\zeta$ is represented as 
\begin{equation*} 
f_{\zeta} (\lambda )
  =\lambda^{\delta} F_{\delta} \left( \lambda +1 +\dfrac{1}{\lambda} \right) 
  =\sum^{2\delta}_{k=0} c_k \lambda^k ,   
\end{equation*} 
where 
\begin{equation*} 
  c_{2\delta} =c_0      =1, \ 
  c_k    =c_{2\delta -k} =\dfrac{c'_k}{2^k} \ 
  (1\leq k\leq \delta -1), \ 
  c_{\delta}  =\dfrac{2c'_{\delta} +1}{2^{\delta}} , 
  c'_k , c'_{\delta} \in \mbox{\boldmath{$Z$}} . 
\end{equation*} 
For $\zeta =\exp (\sqrt{-1} \psi \pi )$ 
with $\psi \in \mbox{\boldmath{$R$}}$, 
suppose that $\psi$ is a rational number. 
Then there exists a positive integer $p$ 
such that $f(\lambda ):=\lambda^{2p} -1$ satisfies $f(\zeta )=0$. 
Therefore we have $f(\lambda)=g(\lambda )f_{\zeta} (\lambda )$ 
and $p>\delta$. 
Then based on the properties of cyclotomic polynomials, 
we have a contradiction and therefore 
we obtain Lemma~\ref{lem:prime}. 
On the other hand, 
we can obtain Lemma~\ref{lem:prime} by direct computations as follows. 
We represent $g$ as $g(\lambda )=\sum^{2p-2\delta}_{k=0} A_k \lambda^k$, 
where $A_0 =-1$, $A_k \in \mbox{\boldmath{$Q$}}$. 
We represent $A_k$ as $A_k =B_k /2^k$ with $B_k \in \mbox{\boldmath{$Z$}}$. 
Suppose $2p-2\delta \geq \delta$. 
If $c'_k$ ($1\leq k\leq \delta -1$) are even, 
then $B_{\delta}$ is odd, and then we see by induction 
that $B_{\delta l}$ is odd for $l\geq 1$. 
This contradicts $A_{2p-2\delta -k} =-A_k$. 
Suppose that $c'_k$ is odd for an integer $k\in \{ 1, \dots , \delta -1\}$ 
and that $c'_l$ ($1\leq l\leq k-1$) are even. 
Then $B_k$ is odd. 
Then there exists $l\in \{ 1, \dots , \delta \}$ 
such that $B_{k+l}$ is odd. 
We see by induction that if $B_{k'}$ is odd, 
then there exists $l\in \{ 1, \dots , \delta \}$ 
such that $B_{k'+l}$ is odd. 
This contradicts $A_{2p-2\delta -k} =-A_k$. 
Even if we suppose $2p-2\delta <\delta$, 
we have a similar contradiction. 
Hence we have concluded that $\psi$ is an irrational number, 
and therefore we obtain Lemma~\ref{lem:prime} again. 
\hfill 
$\square$ 

\vspace{3mm} 

We will prove Proposition~\ref{pro:inftyinfty}. 
By Lemma~\ref{lem:evls}, 
the two elements $C'_{\varepsilon , 1} C'_{\varepsilon , 2}$, 
                 $C'_{\varepsilon , 2} C'_{\varepsilon , 1}$ 
have the same order. 
By Lemmas~\ref{lem:case(ii)infty}, 
          \ref{lem:case(iii)infty}, 
          \ref{lem:case(iv(a))infty}, 
          \ref{lem:case(iv(b))infty}, 
          \ref{lem:case(iv(c))infty} 
and       \ref{lem:prime}, 
the above two elements have infinite order. 
Therefore by Lemma~\ref{lem:1evcts}, 
they are noncommutative. 
Hence we obtain Proposition~\ref{pro:inftyinfty}. 

\begin{remark}\label{rem:case(i)infty}  
Suppose that $n>2$ is an integer and not equal to $4$, 
and set $(\theta'_{\varepsilon , 1} , 
          \theta'_{\varepsilon , 2} , 
          \phi_{\varepsilon} ) 
        =(\pi /2, 2\pi /n , \pi /2)$. 
Then it can be shown that 
the topological holonomy group $G_{\hat{\nabla} , \varepsilon}$ 
of $\hat{\nabla}$ 
in $\bigwedge^2_{\varepsilon}\!E$ at ${\rm pr}(0,0)$ generated by 
$C'_{\varepsilon , 1}$, $C'_{\varepsilon , 2}$ 
has a noncommutative pair of two infinite order elements. 
If $n>2$ is a prime number, 
then this is already shown by Lemma~\ref{lem:prime}. 
In the case where $n>2$ is not prime, 
we observe it as follows. 
\begin{itemize}
\item[{\rm (i)}]{Suppose that $n$ is represented as $n=2^m$ 
with $m\geq 3$. 
We set 
\begin{equation*} 
\hat{C}_1 :=C'_{\varepsilon , 1} (C'_{\varepsilon , 2} )^{2^{m-3}} , \quad 
\hat{C}_2 :=(C'_{\varepsilon , 2} )^{2^{m-3}} C'_{\varepsilon , 1} . 
\end{equation*} 
Then by Lemma~\ref{lem:evls}, 
$\hat{C}_1$, $\hat{C}_2$ have the same order. 
Moreover, the minimal polynomial $f_{\zeta}$ of an eigenvalue $\zeta$ 
of $\hat{C}_k$ other than $1$ on $\mbox{\boldmath{$Q$}}$ is 
given as in \eqref{deg4_2}. 
Therefore $\hat{C}_1$, $\hat{C}_2$ have infinite order. 
Therefore, by Lemma~\ref{lem:1evcts}, 
$\hat{C}_1$, $\hat{C}_2$ are noncommutative.} 
\item[{\rm (ii)}]{Suppose that $n$ has a prime factor $p>2$. 
Then we set 
\begin{equation*} 
\hat{C}_1 :=C'_{\varepsilon , 1} (C'_{\varepsilon , 2} )^{\frac{n}{p}} , 
\quad 
\hat{C}_2 :=(C'_{\varepsilon , 2} )^{\frac{n}{p}} C'_{\varepsilon , 1} . 
\end{equation*} 
Then $\hat{C}_1$, $\hat{C}_2$ have the same order. 
Moreover, referring to the proof of Lemma~\ref{lem:prime}, 
we see that $\hat{C}_1$, $\hat{C}_2$ have infinite order. 
Therefore, by Lemma~\ref{lem:1evcts}, 
$\hat{C}_1$, $\hat{C}_2$ are noncommutative.} 
\end{itemize} 
\end{remark} 

\section{The main theorems}\label{sect:mainthms} 

\setcounter{equation}{0} 

Referring to \eqref{trace}, 
for fixed $\theta_1$, $\theta_2 \in (0, 2\pi )$, 
let $f_{\theta_1 , \theta_2} =f_{\theta_1 , \theta_2} (\phi )$ be 
a function of $\phi \in (0, \pi /2]$ defined by 
\begin{equation} 
\begin{split} 
      f_{\theta_1 , \theta_2} (\phi ) 
:= &     \cos \theta_1  \cos \theta_2 
     - 2 \sin \theta_1  \sin \theta_2  \cos   \phi \\ 
   & +(1+\cos \theta_1  \cos \theta_2 )\cos^2 \phi 
     +(  \cos \theta_1 +\cos \theta_2 )\sin^2 \phi . 
\end{split}
\label{trace2}
\end{equation} 
Then $f_{\theta_1 , \theta_2} (\phi )$ is equal to 
the trace ${\rm tr}\,(C_1 C_2 )$ of $C_1 C_2$ 
with $C_1 =C(\theta_1 )$ and 
     $C_2 =U(\phi )C(\theta_2 ){}^t U(\phi )$. 
Since the absolute values of the eigenvalues of each element of $SO(3)$ is 
equal to one, 
by \eqref{trace2}, 
$f_{\theta_1 , \theta_2} (\phi )$ is contained in the interval $[-1, 3]$. 
Let $I_{\theta_1 , \theta_2}$ be 
the image of the function $f_{\theta_1 , \theta_2}$: 
\begin{equation*} 
I_{\theta_1 , \theta_2} 
 :=\{ f_{\theta_1 , \theta_2} (\phi ) \ | \ \phi \in (0, \pi /2]\} \ 
  (\subset [-1, 3]). 
\end{equation*} 
Let $S$ be a countable set of the interval $[-1, 3]$ 
defined by 
\begin{equation*} 
S:=\{ 1+2\cos (q\pi ) \ | \ q\in \mbox{\boldmath{$Q$}} \} . 
\end{equation*} 
If $(\theta_1 , \theta_2 )\not= (\pi , \pi )$, 
then for each $r\in I_{\theta_1 , \theta_2} \setminus S$ 
and each number $\phi \in (0, \pi /2]$ 
satisfying $f_{\theta_1 , \theta_2} (\phi )=r$, 
$\hat{C}_1 :=C_1 C_2$, 
$\hat{C}_2 :=C_2 C_1$ 
are noncommutative and satisfy ${\rm ord}\,(\hat{C}_k )=\infty$ 
for $k=1, 2$. 
In particular, we obtain 

\begin{theorem}\label{thm:inftyinfty} 
For $k=1, 2$, 
set $\theta_{\varepsilon , k} =q_k \pi$ 
with $(q_1 , q_2 ) 
      \in (\mbox{\boldmath{$Q$}} \times \mbox{\boldmath{$Q$}} ) 
           \setminus 
          (\mbox{\boldmath{$Z$}} \times \mbox{\boldmath{$Z$}} )$ 
and let $\phi_{\varepsilon} \in (0, \pi /2]$ 
satisfy $f_{\theta_{\varepsilon , 1} , \theta_{\varepsilon , 2}} 
             (\phi_{\varepsilon} )\in I_{\theta_1 , \theta_2} \setminus S$. 
Then the topological holonomy group $G_{\hat{\nabla} , \varepsilon}$ is 
generated by 
finite order elements $C_{\varepsilon , 1}$, $C_{\varepsilon , 2}$ and 
equipped with 
a noncommutative pair $(C_{\varepsilon , 1} C_{\varepsilon , 2} , 
                        C_{\varepsilon , 2} C_{\varepsilon , 1} )$ 
of two infinite order elements. 
\end{theorem} 

We will prove 

\begin{theorem}\label{thm:main} 
Let $\nabla$ be an $h$-connection of $E$ 
such that the topological holonomy group $G_{\hat{\nabla} , \varepsilon}$ 
has a noncommutative pair $(\hat{C}_1 , \hat{C}_2 )$ 
with $\infty ={\rm ord}\,(\hat{C}_1 ) 
      \geq    {\rm ord}\,(\hat{C}_2 ) \geq 3$. 
Then for each element $\omega$ of $\hat{E}_{\varepsilon , {\rm pr}(0,0)}$, 
$X(\omega )$ as in \eqref{ch} is 
a dense subset of $\hat{E}_{\varepsilon , {\rm pr}(0,0)}$. 
\end{theorem} 

\vspace{3mm} 

\par\noindent 
\textit{Proof} \ 
By ${\rm ord}\,(\hat{C}_1 ) =\infty$, 
an eigenvalue $\zeta_1$ of $\hat{C}_1$ other than $1$ is 
represented as $\zeta_1 =\exp (\sqrt{-1} \psi_1 \pi )$ 
for a real, irrational number $\psi_1$. 
Let $p_0$ be a point of $S^2$ and set 
\begin{equation*} 
W_1 :=\{ \hat{C}^l_1 p_0 \ | \ 
                 l\in \mbox{\boldmath{$N$}} \cup \{ 0\} \} . 
\end{equation*} 
If $W_1 \not= \{ p_0 \}$, 
then $W_1$ is a dense subset of a circle in $S^2$. 
We set 
\begin{equation*} 
W'_1 :=\{ \hat{C}^m_2 p \ | \ 
                  m\in \mbox{\boldmath{$N$}} \cup \{ 0\} , \ p\in W_1 \} . 
\end{equation*} 
If $W_1 = \{ p_0 \}$, 
then noticing that $\hat{C}_1$, $\hat{C}_2$ are noncommutative, 
we see that $W'_1$ is a subset of a circle in $S^2$ 
which has another point than $p_0$. 
In addition, 
noticing ${\rm ord}\,(\hat{C}_2 ) \geq 3$, 
we see that $W'_1$ has a point of $S^2$ which coincides with 
neither $p_0$ nor $-p_0$. 
Inductively, we set 
\begin{equation*} 
\begin{split} 
W_i  & :=\{ \hat{C}^l_1 p \ | \ l\in \mbox{\boldmath{$N$}} \cup \{ 0\} , \ 
                        p\in W'_{i-1} \} , \\ 
W'_i & :=\{ \hat{C}^m_2 p \ | \ m\in \mbox{\boldmath{$N$}} \cup \{ 0\} , \ 
                        p\in W_i \} 
\end{split} 
\end{equation*} 
for $i\geq 2$. 
Then $W_2$ contains at least one dense subset of a circle in $S^2$. 
In addition, $W'_2$ contains dense subsets of two circles in $S^2$ 
such that the planes containing the circles are not parallel to 
each other. 
Therefore there exists a positive integer $i_0$ 
such that $W_{i_0}$ is a dense subset of $S^2$. 
This means that $X(p_0 )$ as in \eqref{X(p)} is dense in $S^2$. 
Hence we have proved Theorem~\ref{thm:main}. 
\hfill 
$\square$ 

\vspace{3mm} 

Based on Theorem~\ref{thm:main}, we will prove 

\begin{theorem}\label{thm:main2} 
Let $\nabla$ be an $h$-connection of $E$ as in Theorem~\ref{thm:main}. 
Then $G_{\hat{\nabla} , \varepsilon}$ is a dense subset of $SO(3)$. 
\end{theorem} 

Let $\nabla$ be an $h$-connection of $E$ as in Theorem~\ref{thm:main}. 
Let $A=[\mbox{\boldmath{$a$}}_1 \ 
        \mbox{\boldmath{$a$}}_2 \ 
        \mbox{\boldmath{$a$}}_3 ]$ be an arbitrarily given element 
of $SO(3)$. 
In order to prove Theorem~\ref{thm:main2}, 
we will find a sequence $\{ X_n \}_{n\geq 1}$ 
in $G_{\hat{\nabla} , \varepsilon}$ which converges to $A$. 
This is equivalent to the condition that 
if we represent each $X_n$ 
as $X_n =[\mbox{\boldmath{$x$}}_{n1} \ 
          \mbox{\boldmath{$x$}}_{n2} \ 
          \mbox{\boldmath{$x$}}_{n3} ]$, 
then $\{ \mbox{\boldmath{$x$}}_{ni} \}_{n\geq 1}$ converges to 
$\mbox{\boldmath{$a$}}_i$ for any $i\in \{ 1, 2, 3\}$ 
with respect to the standard topology of $\mbox{\boldmath{$R$}}^3$. 

We will prove 

\begin{lemma}\label{lem:eta1} 
Let $\eta_1$ be a unit eigenvector of $\hat{C}_1$ 
corresponding to $\lambda =1$. 
Then there exists a sequence $\{ P_n \}_{n\geq 1}$ 
in $G_{\hat{\nabla} , \varepsilon}$ 
with $P_n =[\mbox{\boldmath{$p$}}_{n1} \ 
            \mbox{\boldmath{$p$}}_{n2} \ 
            \mbox{\boldmath{$p$}}_{n3} ]$ 
satisfying 
\begin{itemize}
\item[{\rm (a)}]{$\{ \mbox{\boldmath{$p$}}_{n1} \}_{n\geq 1}$ converges to 
$\eta_1 ;$} 
\item[{\rm (b)}]{$\{ \mbox{\boldmath{$p$}}_{n2} \}_{n\geq 1}$ converges to 
a unit vector $\eta^{\perp}_1$ orthogonal to $\eta_1$ 
with respect to the standard inner product in $\mbox{\boldmath{$R$}}^3 ;$} 
\item[{\rm (c)}]{$\{ \mbox{\boldmath{$p$}}_{n3} \}_{n\geq 1}$ converges to 
$\eta_1 \times \eta^{\perp}_1$ 
for the standard vector product $\times$ in $\mbox{\boldmath{$R$}}^3$.} 
\end{itemize} 
\end{lemma} 

\vspace{3mm} 

\par\noindent 
\textit{Proof} \ 
Setting $p_0 :=p_1 ={}^t [1 \ 0 \ 0]$ in the proof of Theorem~\ref{thm:main}, 
we see that there exists a sequence $\{ P'_n \}_{n\geq 1}$ 
in $G_{\hat{\nabla} , \varepsilon}$ 
with $P'_n =[\mbox{\boldmath{$p$}}'_{n1} \ 
             \mbox{\boldmath{$p$}}'_{n2} \ 
             \mbox{\boldmath{$p$}}'_{n3} ]$ 
such that $\{ \mbox{\boldmath{$p$}}'_{n1} \}_{n\geq 1}$ converges to 
$\eta_1$. 
Then we can find 
a subsequence $\{ \mbox{\boldmath{$p$}}'_{n(j)2} \}_{j\geq 1}$ 
of $\{ \mbox{\boldmath{$p$}}'_{n2} \}_{n\geq 1}$ which converges to 
a unit vector $\eta^{\perp}_1$ orthogonal to $\eta_1$. 
Since $\mbox{\boldmath{$p$}}'_{n(j)3} 
      =\mbox{\boldmath{$p$}}'_{n(j)1} \times 
       \mbox{\boldmath{$p$}}'_{n(j)2}$, 
$\{ \mbox{\boldmath{$p$}}'_{n(j)3} \}_{j\geq 1}$ converges to 
$\eta_1 \times \eta^{\perp}_1$. 
Therefore it is observed that 
the subsequence $\{ P'_{n(j)} \}_{j\geq 1}$ of $\{ P'_n \}_{n\geq 1}$ is 
suitable for 
$\{ P_n \}_{n\geq 1}$ in Lemma~\ref{lem:eta1}, 
and hence we have proved Lemma~\ref{lem:eta1}. 
\hfill 
$\square$ 

\vspace{3mm} 

For $\eta_1$, $\eta^{\perp}_1$ as in Lemma~\ref{lem:eta1}, 
we set $B:=[\eta_1 \ \eta^{\perp}_1 \ \eta_1 \times \eta^{\perp}_1 ]$. 
Then $B$ is an element of $SO(3)$ and 
the sequence $\{ P_n \}_{n\geq 1}$ in Lemma~\ref{lem:eta1} converges to $B$. 

We will prove 

\begin{lemma}\label{lem:eta1a1} 
Let $\eta_1$, $\eta^{\perp}_1$ be as in Lemma~\ref{lem:eta1}. 
Then there exists a sequence $\{ Q_n \}_{n\geq 1}$ 
in $G_{\hat{\nabla} , \varepsilon}$ 
such that $\{ Q_n B\}_{n\geq 1}$ converges to $A$, that is, 
$\{ Q_n \}_{n\geq 1}$ satisfies 
\begin{itemize}
\item[{\rm (a)}]{$\{ Q_n \eta_1 \}_{n\geq 1}$ converges to 
$\mbox{\boldmath{$a$}}_1 ;$} 
\item[{\rm (b)}]{$\{ Q_n \eta^{\perp}_1 \}_{n\geq 1}$ converges to 
$\mbox{\boldmath{$a$}}_2 ;$} 
\item[{\rm (c)}]{$\{ Q_n (\eta_1 \times \eta^{\perp}_1 )\}_{n\geq 1}$ 
converges to $\mbox{\boldmath{$a$}}_3$.} 
\end{itemize} 
\end{lemma} 

\vspace{3mm} 

\par\noindent 
\textit{Proof} \ 
As in the proof of Theorem~\ref{thm:main}, we suppose that 
an eigenvalue $\zeta_1$ of $\hat{C}_1$ other than $1$ is represented as 
$\zeta_1 =\exp (\sqrt{-1} \psi_1 \pi )$ 
for a real, irrational number $\psi_1$. 
Noticing that $\hat{C}_1$ fixes $\eta_1$ 
and that $\hat{C}_1$ induces the rotation with angle $\psi_1 \pi$ 
around the straight line in $\mbox{\boldmath{$R$}}^3$ 
through the origin and given by $\eta_1$, 
we see that for $\eta^{\perp}_1$ as in Lemma~\ref{lem:eta1}, 
$\{ \hat{C}^m_1 \eta^{\perp}_1 \}_{m\geq 1}$ is a dense subset of 
the great circle which is the intersection of the unit sphere 
in $\mbox{\boldmath{$R$}}^3$ centered at the origin and 
the plane in $\mbox{\boldmath{$R$}}^3$ 
through the origin and orthogonal to $\eta_1$. 
Setting $p_0 :=\eta_1$ in the proof of Theorem~\ref{thm:main}, 
we see that there exists a sequence $\{ Q'_l \}_{l\geq 1}$ 
in $G_{\hat{\nabla} , \varepsilon}$ 
such that $\{ Q'_l \eta_1 \}_{l\geq 1}$ converges to 
$\mbox{\boldmath{$a$}}_1$.  
Since $\mbox{\boldmath{$a$}}_2$ is orthogonal to $\mbox{\boldmath{$a$}}_1$, 
the closure of $\{ Q'_l \hat{C}^m_1 \eta^{\perp}_1 \}_{l, m\geq 1}$ contains 
$\mbox{\boldmath{$a$}}_2$. 
Therefore we can choose a subsequence $\{ Q_n \}_{n\geq 1}$ 
of $\{ Q'_l \hat{C}^m_1 \}_{l, m\geq 1}$ 
so that $\{ Q_n \eta^{\perp}_1 \}_{n\geq 1}$ converges to 
$\mbox{\boldmath{$a$}}_2$. 
Noticing that $\hat{C}_1$ fixes $\eta_1$, 
we see that $\{ Q_n \eta_1 \}_{n\geq 1}$ converges to 
$\mbox{\boldmath{$a$}}_1$. 
Therefore $\{ Q_n (\eta_1 \times \eta^{\perp}_1 )\}_{n\geq 1}$ 
converges to $\mbox{\boldmath{$a$}}_3$. 
Hence we have concluded that $\{ Q_n \}_{n\geq 1}$ satisfies (a), (b), (c) 
in Lemma~\ref{lem:eta1a1} and we obtain Lemma~\ref{lem:eta1a1}. 
\hfill 
$\square$ 

\vspace{3mm} 

We will prove Theorem~\ref{thm:main2}. 
We set $X_n :=Q_n P_n$. 
Then $\{ X_n \}_{n\geq 1}$ is a sequence 
in $G_{\hat{\nabla} , \varepsilon}$. 
Noticing $X_n =(Q_n B)B^{-1} P_n$, 
we see that $\{ X_n \}_{n\geq 1}$ converges to $AB^{-1} B=A$. 
Hence we obtain Theorem~\ref{thm:main2}. 

We can prove the following classification theorem 
of the topological holonomy groups 
in $\bigwedge^2_{\varepsilon}\!E$. 

\begin{theorem}\label{thm:classification} 
The topological holonomy group $G_{\hat{\nabla} , \varepsilon}$ is 
given by one of the following\/$:$ 
\begin{itemize} 
\item[{\rm (a)}]{a finite group,} 
\item[{\rm (b)}]{a commutative infinite group 
which is generated by one or two elements, 
and dense in a subgroup of $SO(3)$ isomorphic to $SO(2)$,} 
\item[{\rm (c)}]{a noncommutative infinite group generated by 
two elements of order $2$, $\infty$ 
such that these rotation axes are orthogonal,} 
\item[{\rm (d)}]{a noncommutative infinite group 
which is dense in $SO(3)$.} 
\end{itemize} 
\end{theorem} 

\vspace{3mm} 

\par\noindent 
\textit{Proof} \ 
In order to prove Theorem~\ref{thm:classification}, 
we can suppose that $G_{\hat{\nabla} , \varepsilon}$ is 
neither finite nor commutative. 
If both of the generators $C_{\varepsilon , 1}$, $C_{\varepsilon , 2}$ 
of $G_{\hat{\nabla} , \varepsilon}$ have order two, 
then by Proposition~\ref{pro:pipi}, 
$G_{\hat{\nabla} , \varepsilon}$ is a group 
as in (c) of Theorem~\ref{thm:classification}. 
Suppose that one of $C_{\varepsilon , 1}$, $C_{\varepsilon , 2}$ is 
of order more than two. 
If the order is $\infty$, 
then by Propositions~\ref{pro:pipi/2}, \ref{pro:2irr} 
or Theorem~\ref{thm:main2}, 
$G_{\hat{\nabla} , \varepsilon}$ is a group 
as in (c) or (d) of Theorem~\ref{thm:classification}. 
Suppose that the order is more than two and finite. 
By Selberg's lemma, 
$G_{\hat{\nabla} , \varepsilon}$ has 
an infinite order element $\hat{C}$. 
If the rotation axis of $\hat{C}$ coincides with 
that of either $C_{\varepsilon , 1}$ or $C_{\varepsilon , 2}$, 
then by Proposition~\ref{pro:2irr} or Theorem~\ref{thm:main2}, 
$G_{\hat{\nabla} , \varepsilon}$ is a group 
as in (d) of Theorem~\ref{thm:classification}. 
If the rotation axis of $\hat{C}$ does not coincide with 
the rotation axis of any of $C_{\varepsilon , 1}$, $C_{\varepsilon , 2}$, 
then we obtain the same conclusion. 
Hence we have proved Theorem~\ref{thm:classification}. 
\hfill 
$\square$ 

\vspace{3mm} 

\begin{remark}\label{rem:tits}  
According to the Tits alternative (see \cite{tits}), 
a finitely generated subgroup of $GL(n, \mbox{\boldmath{$R$}})$ contains 
either a solvable subgroup of finite index 
or     a non-abelian free subgroup. 
We can refer to \cite{BG}, \cite{BG2}, \cite{kuranishi2} 
for the Tits alternative. 
Since the topological holonomy group $G_{\hat{\nabla} , \varepsilon}$ is 
finitely generated and contained in $SO(3)$, 
the Tits alternative is valid for $G_{\hat{\nabla} , \varepsilon}$. 
Then $G_{\hat{\nabla} , \varepsilon}$ 
as in (a), (b) or (c) of Theorem~\ref{thm:classification} is 
in the former type 
(notice that the alternating group $\mathcal{A}_4$ of degree $4$ is 
solvable), 
while $G_{\hat{\nabla} , \varepsilon}$ as in (d) is in the latter type. 
In the case of (d), 
$G_{\hat{\nabla} , \varepsilon}$ is not necessarily isomorphic to 
any free group, 
because there exist topological holonomy groups 
as in (d) generated by finite order elements, 
which are given 
in Section~\ref{sect:ord>2thg} and Theorem~\ref{thm:inftyinfty}. 
Notice that two elements 
$\hat{C}_1 :=C'_{\varepsilon , 1} C'_{\varepsilon , 2}$, 
$\hat{C}_2 :=C'_{\varepsilon , 2} C'_{\varepsilon , 1}$ 
as in Proposition~\ref{pro:inftyinfty} do not necessarily generate 
a non-abelian free subgroup. 
For example, 
if $(\theta'_{\varepsilon , 1} , 
     \theta'_{\varepsilon , 2} , 
     \phi_{\varepsilon} ) 
   =(\pi /3, \pi , \cos^{-1}\!\sqrt{1/3} )$ 
(see Lemma~\ref{lem:case(ii)infty}), 
then ${\rm ord}\,(C'_{\varepsilon , 2})=2$,    
and therefore we have $\hat{C}_1 \hat{C}_2 =(C'_{\varepsilon , 1} )^2$    
and ${\rm ord}\,(\hat{C}_1 \hat{C}_2 )=3$. 
\end{remark} 

Related to (d) of Theorem~\ref{thm:classification}, 
the following problem seems to be interesting, and not easy: 

\begin{problem}
Classify the topological holonomy groups 
as in {\rm (d)} of Theorem~\ref{thm:classification}, 
up to isomorphisms. 
\end{problem} 

Referring to the proof of Theorem~\ref{thm:main2}, we will prove 

\begin{theorem}\label{thm:main3} 
Let $\nabla$ be an $h$-connection of $E$ 
such that for each $\varepsilon \in \{ +, -\}$, 
the generators $C_{\varepsilon , 1}$, $C_{\varepsilon , 2}$ 
of the topological holonomy group $G_{\hat{\nabla} , \varepsilon}$ 
satisfy 
\begin{equation*} 
\theta_{+, 1} =   q_1 \pi , \ 
\theta_{+, 2} =\psi_2 \pi , \ 
\theta_{-, 1} =\psi_1 \pi , \ 
\theta_{-, 2} =   q_2 \pi , \ 
\phi_{\pm} \in (0, \pi /2] 
\end{equation*} 
for $q_1$, $q_2 \in \mbox{\boldmath{$Q$}} \setminus \mbox{\boldmath{$Z$}}$, 
$\psi_1$, $\psi_2 \in \mbox{\boldmath{$R$}} \setminus \mbox{\boldmath{$Q$}}$. 
Then the topological holonomy group $G_{\nabla}$ is dense in $SO(4)$. 
\end{theorem} 

In order to prove Theorem~\ref{thm:main3}, 
we need to show that the subgroup $G_{\hat{\nabla} , +, -}$ 
of $SO(3)\times SO(3)$ 
generated by $(C_{+ , 1} , C_{-, 1} )$, $(C_{+ , 2} , C_{-, 2} )$ 
is dense in $SO(3)\times SO(3)$. 
In order to prove this, we need a lemma. 
Referring to the proof of Theorem~\ref{thm:main}, we will prove 

\begin{lemma}\label{lem:main} 
Let $p_{+0}$, $p_{-0}$ be arbitrarily given points of $S^2$. 
Then 
\begin{equation*} 
\{ (X_+ p_{+0} , X_- p_{-0} ) \ | \ 
   (X_+ , X_- )\in G_{\hat{\nabla} , +, -} \} 
\end{equation*} 
is a dense subset of $S^2 \times S^2$. 
\end{lemma} 

\vspace{3mm} 

\par\noindent 
\textit{Proof} \ 
Let $p_+$, $p_-$ be points of $S^2$. 
Let $O_+$, $O_-$ be neighborhoods of $p_+$, $p_-$ in $S^2$ respectively. 
Referring to the proof of Theorem~\ref{thm:main}, 
we see that there exists an element $(Y_+ , Y_- )$ 
of $G_{\hat{\nabla} , +, -}$ satisfying $Y_+ p_{+0} \in O_+$. 
We set $p'_+ :=Y_+ p_{+0}$, $p'_- :=Y_- p_{-0}$. 
Let $N_+$ be a positive integer such that $N_+ q_1$ is an even integer. 
Then $C^{N_+ l}_{+, 1} =E_3$ 
for any $l\in \mbox{\boldmath{$N$}} \cup \{ 0\}$, 
and $\{ C^{N_+ l}_{-, 1} p'_-  \ | \ 
               l\in \mbox{\boldmath{$N$}} \cup \{ 0\} \}$ 
is either $\{ p'_- \}$ or a dense subset of a circle of $S^2$. 
If $\{ C^{N_+ l}_{-, 1} p'_-  \ | \ 
              l\in \mbox{\boldmath{$N$}} \cup \{ 0\} \} 
   =\{ p'_- \}$, 
then there exists a positive integer $n'$ satisfying 
\begin{itemize} 
\item[{\rm (i)}]{$C^{n'}_{-, 2} p'_- \not= \pm p'_-$,} 
\item[{\rm (ii)}]{$C^{n'}_{+, 2} p'_+ \in O_+$.}  
\end{itemize} 
Therefore in the same way as in the proof of Theorem~\ref{thm:main}, 
we see that there exists an element $(Z_+ , Z_- )$ 
of $G_{\hat{\nabla} , +, -}$ 
satisfying $Z_{\pm}  p'_{\pm} \in O_{\pm}$. 
Hence we have proved Lemma~\ref{lem:main}. 
\hfill 
$\square$ 

\vspace{3mm} 

We will prove Theorem~\ref{thm:main3}. 
By Lemma~\ref{lem:main}, 
there exists a sequence $\{ (P_{+, n} , P_{-, n} )\}_{n\geq 1}$ 
in $G_{\hat{\nabla} , +, -}$ 
with $P_{\varepsilon , n} 
       =[\mbox{\boldmath{$p$}}_{\varepsilon , n1} \ 
         \mbox{\boldmath{$p$}}_{\varepsilon , n2} \ 
         \mbox{\boldmath{$p$}}_{\varepsilon , n3} ]$ satisfying 
\begin{itemize}
\item[{\rm (i)}]{$\{ \mbox{\boldmath{$p$}}_{+, n1} \}_{n\geq 1}$ 
converges to 
an eigenvector of $C_{+ , 2}$ corresponding to $\lambda =1$;} 
\item[{\rm (ii)}]{$\{ \mbox{\boldmath{$p$}}_{-, n1} \}_{n\geq 1}$ 
converges to 
an eigenvector of $C_{- , 1}$ corresponding to $\lambda =1$;} 
\item[{\rm (iii)}]{$\{ \mbox{\boldmath{$p$}}_{\pm , n2} \}_{n\geq 1}$, 
$\{ \mbox{\boldmath{$p$}}_{\pm , n3} \}_{n\geq 1}$ converge.} 
\end{itemize} 
We set $B_{\varepsilon} :=\lim_{n\to \infty} P_{\varepsilon , n}$ 
for $\varepsilon =+, -$. 
Let $N_-$ be a positive integer such that $N_- q_2$ is an even integer. 
Then $C^{N_- l}_{-, 2} =E_3$ 
for any $l\in \mbox{\boldmath{$N$}} \cup \{ 0\}$. 
Let $(A_+ , A_- )$ be an arbitrarily given element of 
$SO(3)\times SO(3)$. 
Then using Lemma~\ref{lem:main} and 
\begin{equation*} 
\{ (C^{N_+ l}_{+, 1} , C^{N_+ l}_{-, 1} ) \ | \ 
  l\in \mbox{\boldmath{$N$}} \cup \{ 0\} \} , \quad  
\{ (C^{N_- l}_{+, 2} , C^{N_- l}_{-, 2} ) \ | \ 
  l\in \mbox{\boldmath{$N$}} \cup \{ 0\} \} , 
\end{equation*} 
we find a sequence $\{ (Q_{+, n} , Q_{-, n} )\}_{n\geq 1}$ 
in $G_{\hat{\nabla} , +, -}$ such that 
$\{ Q_{\pm , n} B_{\pm} \}$ converge to $A_{\pm}$ respectively. 
Therefore $\{ (Q_{+, n} P_{+, n} , Q_{-, n} P_{-, n} )\}_{n\geq 1}$ 
converges to $(A_+ , A_- )$ and 
we observe that $G_{\hat{\nabla} , +, -}$ is dense in $SO(3)\times SO(3)$. 
Then noticing the double covering $SO(4)\longrightarrow SO(3)\times SO(3)$, 
we find an element $A$ of $SO(4)$ and a neighborhood $O$ of $A$ in $SO(4)$ 
such that $G_{\nabla} \cap O$ is dense in $O$. 
This means that the identity element $E_4$ of $SO(4)$ has 
a similar neighborhood in $SO(4)$. 
The set of elements of $SO(4)$ with such neighborhoods is 
non-empty, open and closed in $SO(4)$. 
Since $SO(4)$ is connected, this set coincides with $SO(4)$. 
Hence we have proved Theorem~\ref{thm:main3}. 

Since we already obtain Theorem~\ref{thm:classification}, 
it is natural to have interest in the following problem: 

\begin{problem}\label{prob:classification}
Classify the topological holonomy groups in $SO(4)$, 
up to isomorphisms. 
\end{problem} 

\section{Hermitian vector bundles of complex rank 2 
over tori}\label{sect:HvboverT2} 

\setcounter{equation}{0} 

Let $E$ be a complex vector bundle over $T^2$ of rank $2$. 
Let $h$ be a Hermitian metric of $E$. 
Let $\nabla$ be an $h$-connection of $E$. 
Let $\Tilde{v}_1$, $\Tilde{v}_2$ form 
an orthonormal basis of ${\rm pr}^*\!E_{(0, 0)}$ with respect to $h$. 
For $i=1, 2$, 
let $\Tilde{\xi}_{i, x}$, $\Tilde{\xi}_{i, y}$ be 
parallel sections with respect to $\nabla$ 
of the restrictions of ${\rm pr}^*\!E$ on $l_x$, $l_y$ respectively 
satisfying $\Tilde{\xi}_{i, x} (0, 0) 
           =\Tilde{\xi}_{i, y} (0, 0) 
           =\Tilde{v}_i$. 
There exist elements $B_1$, $B_2 \in U(2)$ satisfying 
\begin{equation} 
\begin{split} 
   (\Tilde{\xi}_{1, x} (a+2\pi , 0)  \ 
    \Tilde{\xi}_{2, x} (a+2\pi , 0)) 
& =(\Tilde{\xi}_{1, x} (a,       0)  \ 
    \Tilde{\xi}_{2, x} (a,       0))B_1 , \\ 
   (\Tilde{\xi}_{1, y} (0, b+2\pi )  \ 
    \Tilde{\xi}_{2, y} (0, b+2\pi )) 
& =(\Tilde{\xi}_{1, y} (0, b      )  \ 
    \Tilde{\xi}_{2, y} (0, b      ))B_2 
\end{split} 
\label{AxAyC} 
\end{equation} 
for any $a$, $b\in \mbox{\boldmath{$R$}}$. 
The \textit{topological holonomy group\/} of $\nabla$ 
at ${\rm pr} (0, 0)$ is the subgroup $G_{\nabla}$ of $U(2)$ 
generated by $B_1$, $B_2$, which is uniquely determined 
up to a conjugate subgroup of $U(2)$. 

We set $\overline{E} :=T^2 \times \mbox{\boldmath{$C$}}^2$. 
This is a product complex bundle over $T^2$. 
The natural Hermitian inner product of $\mbox{\boldmath{$C$}}^2$ gives 
a Hermitian metric $h$ of $\overline{E}$. 
The following proposition is an analogue of 
Proposition~\ref{pro:prod} and similarly obtained. 

\begin{proposition}\label{pro:prod3}
For arbitrarily given two elements $B_1$, $B_2$ of $U(2)$, 
there exists an $h$-connection $\nabla$ of $\overline{E}$ 
with \eqref{AxAyC}. 
\end{proposition} 

Let $I_1$, $I_2$, $I_3$ be elements of the Lie algebra of $SU(2)$ 
defined by 
\begin{equation*} 
I_1 :=\left[ \begin{array}{cc} 
              \sqrt{-1} &         0 \\ 
                     0  & -\sqrt{-1} 
               \end{array} 
      \right] , \ 
I_2 :=\left[ \begin{array}{cc} 
              0 & -1 \\ 
              1 &  0 
               \end{array} 
      \right] , \ 
I_3 :=\left[ \begin{array}{cc} 
                     0  & \sqrt{-1} \\ 
              \sqrt{-1} &        0 
               \end{array} 
      \right] . 
\end{equation*} 
Then the images of $I_1$, $I_2$, $I_3$ 
by the differential of 
the double covering $\Phi :SU(2)\longrightarrow SO(3)$ given by 
\begin{equation} 
\begin{split} 
& \Phi \left( 
       \left[ \begin{array}{cc} 
               b_1 +\sqrt{-1} b_2 & -b_3 +\sqrt{-1} b_4 \\ 
               b_3 +\sqrt{-1} b_4 &  b_1 -\sqrt{-1} b_2 
                \end{array} 
       \right] 
       \right) \\ 
& =
  \left[ \begin{array}{ccc} 
          b^2_1 +b^2_2 - b^2_3 -b^2_4 &  2b_1    b_4   +2b_2    b_3 
                                      & -2b_1    b_3   +2b_2    b_4 \\ 
        -2b_1    b_4   +2b_2    b_3   &   b^2_1 +b^2_3 - b^2_2 -b^2_4 
                                      &  2b_1    b_2   +2b_3    b_4 \\  
         2b_1    b_3   +2b_2    b_4   & -2b_1    b_2   +2b_3    b_4 
                                      &   b^2_1 +b^2_4 - b^2_2 -b^2_3 
            \end{array} 
  \right] 
\end{split} 
\label{Phi} 
\end{equation} 
are $-2J_1$, $-2J_2$, $-2J_3$ respectively 
for $J_1$, $J_2$, $J_3$ as in \eqref{J}, 
where $b_1$, $b_2$, $b_3$, $b_4$ are real numbers satisfying 
$b^2_1 +b^2_2 +b^2_3 +b^2_4 =1$. 
For $\theta \in \mbox{\boldmath{$R$}}$, 
we set $B(\theta ):=\exp (-(\theta /2)I_1 )\in SU(2)$. 
Then $\Phi (B(\theta ))$ coincides with $C(\theta )$ in \eqref{Ctheta}. 
The eigenvalues of $B(\theta )$ are given by 
$e^{\mp \sqrt{-1} \frac{\theta}{2}}$. 
For each element $B$ of $SU(2)$, 
there exist $\theta \in [0, 4\pi )$ and $U\in SU(2)$ 
satisfying $B=UB(\theta ){}^t \overline{U}$. 
We have 
\begin{equation} 
  UB(\theta ){}^t \overline{U} 
    =\exp \left( -\dfrac{\theta}{2} 
               ((|\alpha |^2 -|\beta |^2 )I_1 
                 +2{\rm Im}\,(\alpha \overline{\beta} )I_2 
                 +2{\rm Re}\,(\alpha \overline{\beta} )I_3 ) 
          \right) , 
\label{UBU} 
\end{equation} 
where 
\begin{equation} 
U:=\left[ \begin{array}{cc} 
           \alpha & -\overline{\beta} \\ 
           \beta  &  \overline{\alpha} 
            \end{array} 
   \right] , \quad 
\alpha , \beta \in \mbox{\boldmath{$C$}} , \ 
|\alpha |^2 +|\beta |^2 =1. 
\label{U} 
\end{equation} 
Suppose $B\not= \pm E_2$ for the identity matrix $E_2$. 
Then $\theta \in (0, 4\pi )\setminus \{ 2\pi \}$. 
We set $U=[\mbox{\boldmath{$u$}}_1 \ \mbox{\boldmath{$u$}}_2 ]$. 
Then $\mbox{\boldmath{$u$}}_1$ (respectively, $\mbox{\boldmath{$u$}}_2$) 
is an eigenvector of $B$ corresponding to $e^{-\sqrt{-1} \frac{\theta}{2}}$ 
(respectively, $e^{\sqrt{-1} \frac{\theta}{2}}$). 
For each $i\in \{ 1, 2\}$, 
$\mbox{\boldmath{$u$}}_i$ is considered to be an element of $S^3$. 
By the Hopf fibration $S^3 \longrightarrow \mbox{\boldmath{$C$}}P^1$, 
$\mbox{\boldmath{$u$}}_i$ gives a point of $\mbox{\boldmath{$C$}}P^1$. 
Then the points of $S^2$ 
$\mbox{\boldmath{$u$}}_1$, $\mbox{\boldmath{$u$}}_2$ give through 
a suitable stereographic projection are represented as 
\begin{equation*} 
\pm ((|\alpha |^2 -|\beta |^2 ), \ 2{\rm Im}\,(\alpha \overline{\beta} ), 
                                 \ 2{\rm Re}\,(\alpha \overline{\beta} )) 
\end{equation*} 
(notice \eqref{UBU}). 
Therefore $\mbox{\boldmath{$u$}}_1$, $\mbox{\boldmath{$u$}}_2$ determine 
a one-dimensional subspace of $\mbox{\boldmath{$R$}}^3$, and 
noticing the first column of the matrix in the right side of \eqref{Phi} 
and that $\Phi$ is a homomorphism, we see that 
this is the eigenspace of $\Phi (B)$ corresponding to $\lambda =1$. 

Let $B_1$, $B_2 \in SU(2)$ be as in \eqref{AxAyC}. 
Then choosing $\Tilde{v}_1$, $\Tilde{v}_2$, 
we can suppose $B_1 =B(\theta_1 )$ for $\theta_1 \in [0, 4\pi )$. 
Let $\theta_2 \in [0, 4\pi )$ and $U_2 \in SU(2)$ 
satisfy $B_2 =U_2 B(\theta_2 ){}^t \overline{U}_2$. 
Suppose $B_k \not= \pm E_2$ for $k=1, 2$. 
Let $\phi \in [0, \pi ]$ satisfy 
$\cos \phi =|\alpha_2 |^2 -|\beta_2 |^2$, 
where ${}^t [\alpha_2 \ \beta_2 ]$ is the first column of $U_2$. 
Then choosing $\theta_2 \in (0, 4\pi )\setminus \{ 2\pi \}$, 
we can suppose $\phi \in [0, \pi /2]$. 

The following is an analogue of Proposition~\ref{pro:finite}. 

\begin{proposition}\label{pro:finite2} 
The topological holonomy group $G_{\nabla}$ in $SU(2)$ is finite 
if and only if $\Phi (G_{\nabla} )$ is finite. 
\end{proposition} 

Referring to the last part of the proof of Theorem~\ref{thm:main3}, 
we can prove 

\begin{proposition}\label{pro:dense2} 
The topological holonomy group $G_{\nabla}$ is dense in $SU(2)$ 
if and only if $\Phi (G_{\nabla} )$ is dense in $SO(3)$. 
\end{proposition} 

The following is an analogue of Theorem~\ref{thm:main2}. 

\begin{theorem}\label{thm:main4} 
Suppose that the topological holonomy group $G_{\nabla}$ in $SU(2)$ has 
a noncommutative pair $(\hat{B}_1 , \hat{B}_2 )$ 
with ${\rm ord}\,(\hat{B}_1 )=\infty$ 
and  ${\rm ord}\,(\hat{B}_2 )\not= 2, 4$. 
Then $G_{\nabla}$ is dense in $SU(2)$. 
\end{theorem} 

\vspace{3mm} 

\par\noindent 
\textit{Proof} \ 
Let $(\hat{B}_1 , \hat{B}_2 )$ be as in Theorem~\ref{thm:main4}. 
Then by $\hat{B}^4_k \not= E_2$ for $k=1, 2$, 
$\hat{B}_1$, $\hat{B}_2$ satisfy not only 
the noncommutativity $\hat{B}_1 \hat{B}_2 \not=  \hat{B}_2 \hat{B}_1$ 
but also             $\hat{B}_1 \hat{B}_2 \not= -\hat{B}_2 \hat{B}_1$. 
Therefore $\hat{C}_1 :=\Phi (\hat{B}_1 )$, $\hat{C}_2 :=\Phi (\hat{B}_2 )$ 
form a noncommutative pair 
satisfying $\infty ={\rm ord}\,(\hat{C}_1 ) 
            \geq    {\rm ord}\,(\hat{C}_2 ) \geq 3$. 
Then Theorem~\ref{thm:main2} means 
that $\Phi (G_{\nabla} )$ is dense in $SO(3)$. 
Therefore by Proposition~\ref{pro:dense2}, 
$G_{\nabla}$ is dense in $SU(2)$. 
Hence we have proved Theorem~\ref{thm:main4}. 
\hfill 
$\square$ 

\vspace{3mm} 

Referring to the proof of Theorem~\ref{thm:main3}, 
we will prove 

\begin{theorem}\label{thm:main5} 
Let $E$, $h$, $\nabla$ be as in the beginning of this section. 
Suppose that $G_{\nabla}$ contains two elements $B_1$, $B_2 \in U(2)$ 
satisfying 
\begin{equation*} 
B_1 =B(\psi \pi ), \quad 
B_2 =e^{\sqrt{-1} \gamma} UB(q\pi ){}^t \overline{U} , \quad 
\phi \in (0, \pi /2] 
\end{equation*} 
for $q\in \mbox{\boldmath{$Q$}} \setminus \mbox{\boldmath{$Z$}}$, 
$\psi$, $\gamma \in \mbox{\boldmath{$R$}} \setminus \mbox{\boldmath{$Q$}}$, 
where $U$ is as in \eqref{U} and 
satisfies $\cos \phi =|\alpha |^2 -|\beta |^2$. 
Then the topological holonomy group $G_{\nabla}$ is dense in $U(2)$. 
\end{theorem} 

\vspace{3mm} 

\par\noindent 
\textit{Proof} \ 
Let $N$ be a positive integer such that $Nq/2$ is an even integer. 
Then $B^N_2 =e^{\sqrt{-1} N\gamma} E_2$, 
which implies that $UB(q\pi ){}^t \overline{U}$ is contained in the closure 
of $G_{\nabla}$. 
We see from Theorem~\ref{thm:main4} that 
$B_1$, $B_2 /e^{\sqrt{-1} \gamma}$ generate a dense subgroup in $SU(2)$. 
Since $B^N_2$ generates a dense subgroup 
in the subgroup $\{ e^{\sqrt{-1} t} E_2 \ | \ t\in [0, 2\pi )\}$ of $U(2)$, 
$G_{\nabla}$ is dense in $U(2)$. 
Hence we have proved Theorem~\ref{thm:main5}. 
\hfill 
$\square$ 

\vspace{3mm} 

The following problem is an analogue of Problem~\ref{prob:classification}. 

\begin{problem}
Classify the topological holonomy groups in $U(2)$, 
up to isomorphisms. 
\end{problem} 

\section*{Acknowledgments}

The authors would like to express their cordial gratitude 
to Professors Shin Kato, Masato Mimura and Ikuo Satake 
for valuable discussions. 
Naoya Ando was supported by JSPS KAKENHI Grant Number JP21K03228. 



\vspace{4mm} 

\par\noindent 
\footnotesize{Naoya Ando \\ 
              Faculty of Advanced Science and Technology, 
              Kumamoto University \\ 
              2--39--1 Kurokami, Chuo-ku, Kumamoto 860--8555 Japan} 

\par\noindent  
\footnotesize{E-mail address: andonaoya@kumamoto-u.ac.jp} 

\vspace{4mm} 

\par\noindent 
\footnotesize{Anri Yonezaki \\ 
              Graduate School of Science and Technology, 
              Kumamoto University \\ 
              2--39--1 Kurokami, Chuo-ku, Kumamoto 860--8555 Japan} 

\par\noindent  
\footnotesize{E-mail address: 223d8007@st.kumamoto-u.ac.jp} 


\begin{thebibliography}{99} 

\bibitem{AK}N.\,Ando and T.\,Kihara, 
Horizontality in the twistor spaces associated with vector bundles of rank 4 
on tori, 
J. Geom. {\bf 112} (2021) 19, 26 pages. 

\bibitem{blichfeldt}H.\,F.\,Blichfeldt, 
Finite collineation groups, 
The University of Chicago Press, Chicago, 1917. 

\bibitem{BG}E.\,Breuillard and T.\,Gelander, 
On dense free subgroups of Lie groups, 
J. Algebra {\bf 261} (2003) 448--467. 

\bibitem{BG2}E.\,Breuillard and T.\,Gelander, 
A topological Tits alternative, 
Ann. of Math. (2) {\bf 166} (2007) 427--474. 

\bibitem{bryant}R.\,Bryant, 
Conformal and minimal immersions of compact surfaces into the $4$-sphere, 
J. Differential Geom. {\bf 17} (1982) 455--473. 

\bibitem{cassels}J.\,W.\,S.\,Cassels, 
Local fields, 
Cambridge University Press, Cambridge, 1986. 

\bibitem{ES}J.\,Eells and S.\,Salamon, 
Twistorial construction of harmonic maps of surfaces 
into four-manifolds, 
Annali della Scuola Normale Superiore di Pisa, Classe di Scienze 
{\bf 12} (1985) 589--640. 

\bibitem{friedrich}T.\,Friedrich, 
On surfaces in four-spaces, 
Ann. Glob. Anal. Geom. {\bf 2} (1984) 257--287. 

\bibitem{hitchin}N.\,J.\,Hitchin, 
Polygons and gravitons, 
Math. Proc. Cambridge Philos. Soc. {\bf 85} (1979) 465--476. 

\bibitem{klein}F.\,Klein, 
Vorlesungen \"{u}ber das Ikosaeder und die Aufl\"{o}sung 
der Gleichungen vom f\"{u}nften Grade, 
edited, with an introduction and commentary by P.\,Slodowy, 
Birkh\"{a}user Verlag and B.\,G.\,Teubner, 1993. 

\bibitem{kuranishi2}M.\,Kuranishi, 
On everywhere dense embedding of free groups in Lie groups, 
Nagoya Math. J. {\bf 2} (1951) 63--71. 

\bibitem{nagata}M.\,Nagata, 
Field theory, 
Monographs and textbooks in pure and applied mathematics {\bf 40}, 
Marcel Dekker, Inc, New York and Basel, 1977. 

\bibitem{selberg}A.\,Selberg, 
On discontinuous groups in higher-dimensional symmetric spaces, 
Contributions to function theory, pp. 147--164, 
Tata Institute of Fundamental Research, Bombay, 1960.  

\bibitem{tits}J.\,Tits, 
Free subgroups in linear groups, 
J. Algebra {\bf 20} (1972) 250--270. 

\end{thebibliography}
\end{document}